\documentclass[10pt, twoside]{article}
\usepackage{amsfonts}
\usepackage{amssymb}
\usepackage{enumerate}
\usepackage{epsfig}
\usepackage{caption2}
\usepackage{amsmath}
\usepackage{ntheorem}
\usepackage{alexandrov_ciraolo} 

\renewcommand{\caption}{\fcaption}

\begin{document}
\setlength{\textheight}{7.8truein}  

\runninghead{Wave Propagation in a 3-D Optical Waveguide}{Wave Propagation in a 3-D Optical Waveguide}
\normalsize\textlineskip
\thispagestyle{empty}
\setcounter{page}{1}
\vspace*{0.88truein}

\fpage{1}
\centerline{\bf WAVE PROPAGATION IN A 3-D OPTICAL WAVEGUIDE}
\baselineskip=13pt
\vspace*{0.37truein}

\centerline{\footnotesize OLEG ALEXANDROV}
\vspace*{0.015truein}

\centerline{\footnotesize\it School of Mathematics, University of  Minnesota}
\baselineskip=10pt
\centerline{\footnotesize\it 206 Church Str SE, Minneapolis, MN 55455, USA}
\baselineskip=10pt
\centerline{\footnotesize\it aoleg@math.umn.edu}

\vspace*{10pt}
\centerline{\footnotesize GIULIO CIRAOLO}
\vspace*{0.015truein}

\centerline{\footnotesize\it Dipartimento di Matematica U. Dini, University of Firenze}
\baselineskip=10pt
\centerline{\footnotesize\it viale Morgagni 67/A, 50134 Firenze, Italy}
\baselineskip=10pt
\centerline{\footnotesize\it ciraolo@math.unifi.it}
\vspace*{1.0truein}

\abstracts{In this paper we study the problem of wave propagation in a
3-D optical fiber. The goal is to obtain a solution for the time-harmonic
field caused by a source in a cylindrically symmetric waveguide.  The
geometry of the problem, corresponding to an open waveguide, makes the
problem challenging.  To solve it, we construct a transform theory which
is a nontrivial generalization of a method for solving a 2-D version of
this problem given by Magnanini and Santosa.\cite{MS}
}
{The extension to 3-D is made complicated by the fact that the resulting
eigenvalue problem defining the transform kernel is singular both at the
origin and at infinity.  The singularities require the investigation of
the behavior of the solutions of the eigenvalue problem.  Moreover, the
derivation of the transform formulas needed to solve the wave
propagation problem involves nontrivial calculations.
}
{The paper provides a complete description on how to construct the
solution to the wave propagation problem in a 3-D optical waveguide with
cylindrical symmetry. A follow-up article will study the particular
cases of a step-index fiber and of a coaxial waveguide. In those cases
we will obtain concrete formulas for the field and numerical
examples. }

\vspace*{6pt}
\keywords{Wave propagation, Optical waveguides, Helmholtz equation,
  Green's function, Spectral representation.}
\vspace*{1pt}\textlineskip           
\section{Introduction}\label{sec1}   
\vspace*{-0.5pt}

In this paper, we study the wave propagation in a cylindrical optical
fiber.\footnote{We will use the terms optical waveguide and optical
fiber interchangeably.} As model equation, we use the \emph{Helmholtz
equation}
\begin{equation}\label{1.1}
\Delta u + k^2 n(x, y, z)^2 u = f(x, y, z),\quad (x, y, z) \in\mathbb R^3,
\end{equation}
also called the \emph{time-harmonic wave equation}. The number $k$
 is called the \emph{wavenumber}, and the
function $f$ represents a source of energy. We require that the index of
refraction $n(x,y,z)$ have the form:
$$ n(x,y,z)=\begin{cases} n_{co}(x^2+y^2) , & \textmd{if }
x^2+y^2< R^2 , \\ n_{cl} , & \textmd{if } x^2+y^2 \ge R^2 ,
\end{cases} $$
where $R$ is the radius of the waveguide.

\par The main result of this paper is the construction of a
representation formula for a solution $u$ of \eqref{1.1}
satisfying suitable radiation conditions. Our results
generalize a similar formula obtained by Magnanini and Santosa\cite{MS}
in the two dimensional case.

\par In \cite{MS} it is shown that the energy of the electromagnetic
field is divided into two parts: a part propagates inside the
waveguide as a finite number of distinct guided modes, while the other
part either decays exponentially along the fiber or is radiated
outside.  Our case reveals a new feature: for special choices of the
parameters, new kinds of guided modes appear which, rather than
decaying exponentially outside the fiber, vanish as a power of the
distance from the fiber's axis.

\par As in,\cite{MS} we use the technique of separation of variables.
But there are a number of important differences. We separate
variables in the three cylindrical coordinates $r,$ $\vartheta,$ and
$z.$ The component in $\vartheta$ is simple to solve for, and will
only introduce a sum over $m\in \mathbb Z$. The separation in $z$ can be
solved as in.\cite{MS}  More complicated is the study of the
$r$-coordinate: we will obtain a differential equation in $r$ which is
Bessel's equation, save for a non-constant coefficient.  This equation
will have a singularity at $r=0.$ Also, instead of working as in 2-D
with sine and cosine functions whose zeros are uniformly distributed,
we will work with Bessel's functions whose distribution of zeros is more
complicated. These differences make the method used in \cite{MS}
inapplicable.

Instead, we use the theory of singular self-adjoint eigenvalue
problems for second order differential equations as presented in
\cite{CL} and. \cite{Ti} Applying this theory to our problem and
doing the explicit calculations took some effort, chiefly for the
reason that our self-adjoint eigenvalue problem has (just like in the
2-D case) a coefficient which is, with some restrictions, a general
function, so we cannot obtain solutions for our equation in terms of
concrete functions.

The paper is organized as follows. In section 2 we derive the
second order self-adjoint eigenvalue problem associated with
the Helmholtz equation \eqref{1.1}. In section 3 we will prove a set of
technical lemmas aimed at studying the behavior of the solutions of
this eigenvalue problem. In section 4 we will classify the
solutions of the eigenvalue problem which are ``well-behaved'' (in a
sense to be specified there) as $r\to0$ and $r\to\infty.$ The
motivation is that the electromagnetic field in the fiber will have a
representation in terms of the ``well-behaved'' solutions of this
eigenvalue problem.  In section 5 we summarize the theory of
self-adjoint eigenvalue problems as exposed in \cite{Ti} and.\cite{CL}

The functions defined in section 5 are calculated in section
6.  In section 7 the transform defined in section
5 is computed.  The obtained transform is used in section
8 to find Green's function for the Helmholtz equation
\eqref{1.1}, and in turn, to find the desired electromagnetic
field in the fiber given the source.

\section{The Eigenvalue Problem}\label{sec2}
A typical optical fiber is a cylindrical dielectric waveguide, made of
silica glass or plastic. Its central region is called \emph{core},
surrounded by \emph{cladding}, which has a slightly lower index of
refraction. The cladding is surrounded by a protective \emph{jacket}.
Most of the electromagnetic radiation propagates along the
core.  The electromagnetic field intensity in the cladding decays
exponentially along the radial direction. This is why, the radius of
the cladding, which is typically several times larger than the radius
of the core, can be considered infinite.

\begin{figure}[h]\label{fig1}
\begin{center}
\includegraphics[width=0.8\textwidth]{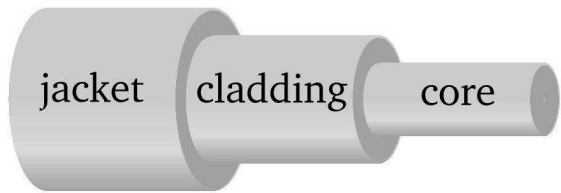}
\caption{An optical fiber.}
\end{center}
\end{figure}

A typical optical fiber is \emph{weakly guided}, which means that the
difference in the indices of refraction of its the core and cladding is
very small. In these conditions the electromagnetic field in the fiber
is essentially \emph{transverse} with each of the transverse components
approximately satisfying the Helmholtz equation \eqref{1.1}.
This is the so-called \emph{weakly guided approximation}. For a
derivation of \eqref{1.1} from Maxwell's equations for a weakly guided
fiber see,\cite{SL} chapters 30 and 32.

Because of the cylindrical geometry it will be convenient to use the
cylindrical coordinate system $(r,\vartheta, z),$  with $z$ being the axial direction.
Then, the index of
refraction will depend on the $r$ variable only, $n=n(r).$ In the new
coordinates  \eqref{1.1} becomes
\begin{equation}\label{2.1}
 \frac{\partial^2 u}{\partial z^2}+\frac{1}{r}\frac{\partial}{\partial r}
 \bigg(r \frac{\partial u}{\partial r}\bigg)+
 \frac{1}{r^2}\frac{\partial^2 u}{\partial \vartheta^2}+k^2 n(r)^2 u=f(r,\vartheta, z).
\end{equation}
This is a linear partial differential equation. The solution of this
equation will be determined as soon as we find its \emph{Green's
function.} In order to obtain the latter, consider the homogeneous
version of \eqref{2.1},
\begin{equation*}
 \frac{\partial^2 u}{\partial z^2}+\frac{1}{r}\frac{\partial}{\partial
 r} \bigg(r \frac{\partial u}{\partial r}\bigg)+
 \frac{1}{r^2}\frac{\partial^2 u}{\partial \vartheta^2}+k^2 n(r)^2
 u=0.
\end{equation*}
Look for a solution in separated variables, $u(r, \vartheta, z)=Z(z)\, \Theta(\vartheta)\,  v(r).$
It is quickly found that we must have
\begin{equation}\label{2.2}
  u(r,\vartheta, z)=e^{i\beta k z}e^{im\vartheta}v(r),
\end{equation}
with $\beta\in\mathbb C,$ $m\in\mathbb Z,$ and $v(r)$ satisfying the
differential equation
\begin{equation*}
  v''+\frac{1}{r}v'+\bigg\{k^2n(r)^2-\beta^{2}-\frac{m^2}{r^2}\bigg\}v=0
\end{equation*}
(the derivative here, and in the rest of this paper will always be in
respect to the variable $r$).

Let $R>0$ be the radius of the fiber core.  Then $n(r)=n_{cl}$ for
$r\ge R.$ Denote
\begin{equation}\label{2.3}
  d^2=n_0^2-n_{cl}^2,\ \  l=k^{2}(n_0^2-\beta^2), \ \  q(r)=k^2[n_0^2-n(r)^{2}].
\end{equation}
Then this equation becomes
\begin{equation}\label{2.4}
  v''+\frac{1}{r}v'+\bigg\{l-q(r)-\frac{m^2}{r^2}\bigg\}v=0.
\end{equation}
We will view \eqref{2.4} as an eigenvalue problem in $l\in\mathbb C.$ The
variable $r$ is in $(0,\infty),$ the number $m$ is an integer, the
function $q$ is bounded, measurable, real-valued and non-negative,
with $q(r)=d^2>0$ for $r\ge R>0$.  It will be convenient to make a
variable change. Denote $w(r)=\sqrt{r}v(r).$ Get the equation
\begin{equation}\label{2.5}
  w''+\bigg\{l - q(r) - \frac{m^2-1/4}{r^2}\bigg\}w=0, \ \ r \in (0, \infty).
\end{equation}
This will be our self-adjoint eigenvalue problem.

\begin{figure}[h]\label{fig2}
\begin{center}
\includegraphics[width=0.8\textwidth]{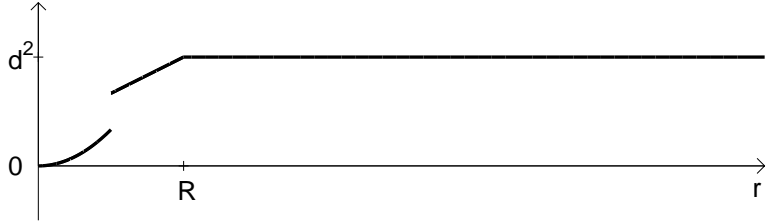}
\caption{The function q(r).}
\end{center}
\end{figure}

\section{A Study of the Solutions of the Eigenvalue Problem}\label{sec3}
Before going further we will need some information about the behavior
of the solutions of the differential equation \eqref{2.5} as functions
of $r$ and $l$ and $m.$ This will be the subject of the next four
lemmas.

\begin{lemma}\label{lm3.1}
There exists a solution $j_m(r, l)$ $(r > 0,\, l\in\mathbb C,\, m\in\mathbb Z)$
of {\rm \eqref{2.5}} such that
\begin{equation}\label{3.1}
  \lim\limits_{r\to0}\frac{j_m(r, l)}{r^{|m|+1/2}}=1, \ \
  \lim\limits_{r\to0}\frac{j_m'(r, l)}{(|m|+1/2)\, r^{|m|-1/2}}=1.
\end{equation}
The functions $j_m(r, l)$ and $j_m'(r, l)$ are analytic in $l$ as $r$
is fixed.  There exists another solution $y_m(r,l)$ of {\rm \eqref{2.5}} such
that
\begin{subequations}
  \begin{equation}\label{3.2a}
    \lim\limits_{r\to0}\frac{y_m(r, l)}{r^{-|m|+1/2}}=1,\ \
    \lim\limits_{r\to0}\frac{y_m'(r, l)}{(-|m|+1/2)\, r^{-|m|-1/2}}=1, \mbox{ if } |m| \ge 1,
  \end{equation}
  and
  \begin{equation}\label{3.2b}
    \lim\limits_{r\to0}\frac{y_m(r, l)}{\sqrt{r}\ln{r}}=1, \ \
    \lim\limits_{r\to0}\frac{y_m'(r, l)}{\ln{r}/(2\sqrt{r})}=1, \mbox{ if } m=0.
  \end{equation}
\end{subequations}
\end{lemma}
\begin{proof}
We will assume $m\ge 0,$ and then set $j_{-m}=j_m$ and $y_{-m}=y_{m}.$
Make the variable change $w=r^{m+1/2}\sigma$ in \eqref{2.5}.
Obtain
$$
  \sigma''+\frac{2m+1}{r}\sigma'+\{l-q(r)\}\,\sigma=0.
$$
Denote $k=2m+1,$ $k\ge 1.$ Multiply this equation by $r^{k}.$ Get
\begin{equation}\label{3.3}
  (r^{k}\sigma')'=r^{k}(q(r)-l)\,\sigma.
\end{equation}
To prove this lemma we need to find two solutions $\sigma(r,l)$ and $\tau(r,l)$
of \eqref{3.3} such that
\begin{equation}\label{3.4}
    \lim\limits_{r\to0}\sigma(r,l)=1, \ \ \lim\limits_{r\to0}\sigma'(r,l)=0,
\end{equation}
\begin{subequations}
\begin{equation}\label{3.5a}
  \lim\limits_{r\to0}r^{2m}\tau(r, l)=1,\ \
  \lim\limits_{r\to0}r^{2m+1}\tau'(r, l)=-2m, \ \ \mbox{ if } m \ge 1,
\end{equation}
and
\begin{equation}\label{3.5b}
  \lim\limits_{r\to0}\frac{\tau(r, l)}{\ln{r}}=1, \ \
  \lim\limits_{r\to0}r\tau'(r, l)=1, \ \ \mbox{ if } m=0.
  \end{equation}
\end{subequations}
Also, $\sigma(r,l)$ must be analytic in $l$ for $r$ fixed.  The idea
to finding $\sigma$ and $\tau$ is to rewrite \eqref{3.3} as an integral
equation.

Let $q_\infty=\sup_{r\in[0,\infty)}q(r).$ If
$\omega:[0,\infty)\to\mathbb C$ is a function, bounded and integrable
on every compact subset of $[0,\infty),$ and $p\ge 0$ is an integer,
then one has
\begin{equation}\label{3.6}
  \bigg|\int\limits_{0}^{r}\!s^{p}\omega(s)\,ds\bigg| \le
  \sup_{t\in[0,r]}|\omega(t)|\int\limits_{0}^{r}\!s^{p}\,ds
  =\sup_{t\in[0,r]}|\omega(t)|\;\frac{r^{p+1}}{p+1}.
\end{equation}
Consider the operator
\begin{equation}\label{3.7}
  T\omega(r)=\int\limits_{0}^{r}\!t^{-k}\int\limits_{0}^{t}\!
  \!s^{k}(q(s)-l)\omega(s)\,ds\,dt,
\end{equation}
defined for complex-valued functions $\omega$ which are bounded and
integrable on every compact subset of $[0,\infty).$ By applying
\eqref{3.6} to the inner-most integral in \eqref{3.7} we deduce
$$
  \sup_{t\in[0,r]}|T\omega(t)|\le (q_\infty+|l|)\sup_{t\in[0,r]}|\omega(t)|\frac{r^{2}}{2\,(k+1)}.
$$
In particular, for any $r>0,$ $T$ is a bounded linear operator from
the space ${\mathcal C}([0,r])$ of continuous, complex-valued
functions defined on $[0, r]$ onto itself.  By using the same
reasoning as above, one can show by induction that for any integer
$n\ge 0$
\begin{equation}\label{3.8}
  |T^n \omega(r)|\le(q_\infty+|l|)^n\sup_{t\in[0, r]}|\omega(t)|
  \frac{r^{2n}}{\{2\cdot 4 \cdots 2n\}\,\{(k+1)\cdot (k+3) \cdots (k+2n-1)\}}.
\end{equation}

Using \eqref{3.6} it is easy to check that if function $\sigma$
satisfies \eqref{3.4}, then \eqref{3.3} is equivalent to
\begin{equation}\label{3.9}
  \sigma'(r,\lambda)=\frac{1}{r^{k}}\int\limits_{0}^{r}\!s^{k}(q(s)-\lambda)\,\sigma(s,\lambda)\,ds,
\end{equation}
which in turn is the same as
\begin{equation}\label{3.10}
  \sigma=1+T\sigma.
\end{equation}

We will find  $\sigma$ by applying the method of successive
approximations. Let $\sigma_0\equiv 1,$ and $\sigma_{n+1}=\sigma_0+T\sigma_n,\,
n\ge0.$ Then,
$$
  \sigma_n=\sigma_0+T\sigma_0+\dots+T^n\sigma_0, \ \ n\ge 0.
$$
By using \eqref{3.8} we obtain that the series
\begin{equation}\label{3.11}
  \sigma=\sigma_0+T\sigma_0+\dots+T^n\sigma_0+\dots, \ \ \sigma_0\equiv 1,
\end{equation}
is uniformly convergent for $r$ and $l$ in compact sets. Since
$\sigma_n$ is continuous in $r$ and analytic in $l,$ the same property
will hold for $\sigma.$ Using the fact that the operator $T$ is
continuous on ${\mathcal C}([0,r])$ for any $r>0$, it follows from the
above series that \eqref{3.10} holds, and therefore \eqref{3.3} holds.

Now prove the existence of $\tau$ satisfying \eqref{3.3}
with the boundary conditions \eqref{3.5a} or \eqref{3.5b}. That will be
done in several steps to be described below.  For the rest of the
lemma, $\sigma=\sigma(r,l)$ will be the solution to \eqref{3.3} found
above. We are not looking for any properties of $\tau$ in respect with
$l.$ Therefore, we will fix $l\in\mathbb C,$ and will consider $\tau$
to be a function of $r$ only. At step 1 we will prove that if $\tau$
is a solution to \eqref{3.3} which is linearly independent of $\sigma,$ then
$\tau$ is unbounded as $r\to0.$ At step 2, we will show that any
solution $\tau$ to \eqref{3.3} has the property
\begin{equation}\label{3.12}
 |\tau(r)|\le A+Br^{-k}, 0< r \le r_0,
\end{equation}
for some $A>0, \, B> 0, r_0> 0.$ At step 3 we will prove that any
solution $\tau$ of \eqref{3.3} linearly independent of $\sigma$
satisfies
\begin{equation}\label{3.13}
  \tau'(r)=O(r^{-k}) \mbox{ as } r\to0.
\end{equation}
That will imply \eqref{3.5a} or \eqref{3.5b} depending on whether $m\ge1$
or $m=0,$ that is, $k\ge 3$ or $k=1.$

{\bf Step 1.} Show that any solution $\tau$ to \eqref{3.3} which is
linearly independent of $\sigma$ is unbounded as $r\to 0.$ Assume that
$\tau$ is such a solution and that it is bounded as $r\to0.$ From
\eqref{3.3} we get
$$
  r^{k}\tau'(r)=c_1+\int\limits_{0}^{r}\!s^k(q(s)-l)\tau(r)\,ds
$$
for some constant $c_1.$ We cannot have $c_1\ne0,$ because then
$\tau'(r)=O(r^{-k})$ with $k\ge 1,$ so $\tau$ cannot be bounded as
$r\to0.$ Then, if $c_1=0$ we can divide by $r^k$ and integrate again,
to obtain
$$
  \tau=c_2+T\tau,
$$
with $T$ the operator defined by \eqref{3.7} and $c_2$ another
constant. Since the solution $\sigma$ satisfies $\sigma=1+T\sigma,$ we
deduce that $\varphi=\tau-c_2\sigma$ will satisfy
$$
  \varphi=T\varphi.
$$
Then $\varphi=T^n\varphi$ for any $n\ge 0.$ By using \eqref{3.8} we get
$$
  |\varphi(r)|\le(q_\infty+|l|)^n\sup_{t\in[0, r]}|\varphi(t)|
  \frac{r^{2n}}{\{2\cdot 4 \cdots 2n\}\,\{(k+1)\cdot (k+3) \cdots (k+2n-1)\}},
$$
for any $r>0$ and $n\ge 0,$ which implies $\varphi=0,$ that is, $\tau=c_2\sigma.$ This
is a contradiction with the assumption that $\tau$ is linearly
independent of $\sigma.$

{\bf Step 2.} Show that for any  $\tau$ solution of \eqref{3.3}
inequality \eqref{3.12} holds.  Let us again write \eqref{3.3} as an
integral equation. This time we cannot integrate from $r=0,$ since we
expect an unbounded solution as $r\to 0$.  Let $r_0>0$ be a fixed number. It is
easy to show that $\tau$ will satisfy \eqref{3.3} if and only if
\begin{equation}\label{3.14}
  \tau=\tau_0+S\tau,
\end{equation}
where
$$
  \tau_0(r)=c_0+d_0\int\limits_{r_0}^{r}\!t^{-k}\,dt,
$$
with $c_0=\tau(r_0),$ $d_0=r_0^{k}\tau'(r_0)$, and the operator $S$ is
defined on functions $\theta$ integrable on compact subsets of $(0,\infty)$
and is given by the formula
$$
  S\theta(r)=\int\limits_{r_0}^{r}\!t^{-k}\int\limits_{r_0}^{t}\!s^k(q(s)-l)\theta(s)\,ds\,dt.
$$

Use once again the method of successive approximations to express $\tau$ as the sum of an infinite series. Let
$c_0\in \mathbb R, \, d_0\in \mathbb R,$ and define
$\tau_0(r)=c_0+d_0\int_{r_0}^{r}\!t^{-k}\,dt,$
$\tau_{n+1}=\tau_0+S\tau_{n}, \, n \ge 0.$ We have
$$
  \tau_n=\tau_0+S\tau_0+\dots+S^n \tau_0, \ \ n\ge 0.
$$

Notice that if $\omega$ is an integrable function on $(0, r_0),$ then for $0< r \le r_0$
\begin{equation}\label{3.15}
  \bigg|\int\limits_{r_0}^{r}\!s^{-k}\omega(s)\,ds\bigg|\le
  \int\limits_{r}^{r_0}\!s^{-k}|\omega(s)|\,ds\le r^{-k}\int\limits_{r}^{r_0}\!|\omega(s)|\,ds.
\end{equation}
By applying \eqref{3.15} we get the estimate
$|\tau_0(r)|\le|c_0|+|d_0|r_0r^{-k}$ for $0<r\le r_0.$ Then,
$$
  \bigg|\int\limits_{r_0}^{t}\!s^k(q(s)-l)\tau_0(s)\,ds\bigg|\le
  \int\limits_{t}^{r_0}\!(q_\infty+|l|)(r_0^k|c_0|+|d_0|r_0)\,ds=C(q_\infty+|l|)(r_0-t),
$$
where $C=r_0^k|c_0|+|d_0|r_0.$ Apply \eqref{3.15} to estimate $S\tau_0.$
$$
  |S\tau_0(r)|\le r^{-k}\int\limits_{r}^{r_0}\!C(q_\infty+|l|)(r_0-t)\,dt =
  C(q_\infty+|l|)r^{-k}\frac{(r_0-r)^2}{2!}.
$$
By repeating the above reasoning, one obtains by induction
$$
  |S^n \tau_0(r)|\le C(q_\infty+|l|)^n r^{-k}\frac{(r_0-r)^{2n}}{(2n)!}, \ \ 0< r\le r_0,\, n\ge 1.
$$
As it was done for $\sigma,$ one can prove that the series
$$
  \tau=\sum\limits_{n=0}^{\infty}S^{n}\tau_0
$$
is uniformly convergent for $r$ in compact subsets of $(0,r_0]$
and that $\tau$ satisfies \eqref{3.14}. One can then estimate $\tau.$
$$
 |\tau(r)|\le \sum\limits_{n=0}^{\infty}|S^{n}\tau_0(r)|\le |c_0|+Cr^{-k}e^{\sqrt{q_\infty+|l|}(r_0-r)}, \ \
  0< r\le r_0,
$$
which implies \eqref{3.12}.

{\bf Step 3.} Show that any solution $\tau$ to \eqref{3.3} which is
linearly independent of $\sigma$ satisfies \eqref{3.13}.  From
\eqref{3.12} we deduce that $r^{k}(q(r)-l)\tau(r)$ is bounded as
$r\to0.$ Integrate \eqref{3.3} from $0$ to $r.$ Get
\begin{equation}\label{3.16}
  r^{k}\tau'(r)=c_3+\int\limits_{0}^{r}\!s^{k}(q(s)-l)\tau(s)\,ds,
\end{equation}
with $c_3$ a constant.

Suppose that $c_3=0.$ Then, from \eqref{3.12} and \eqref{3.16} it follows that
$|r^{k}\tau'(r)|\le Mr$ for $0< r \le r_0$ and some $M>0.$ But
then, $|\tau'(r)|\le M/r^{k-1}$ for $0< r \le r_0.$ Note that $k$ is
an odd number, since $k=2m+1$ with $m\ge 0$ an integer. In particular
$k\ne 2,$ or $k-1\ne 1.$ Since
$\tau(r)=\tau(r_0)+\int_{r_0}^{r}\!\tau'(s)\,ds,$ we obtain
$$
  |\tau(r)|\le A_1+\frac{B_1}{r^{k-2}},
$$
for some $A_1>0,\, B_1>0.$ This is the same as \eqref{3.12} but with
$k-2$ instead of $k$. By repeating several times the same reasoning
starting with \eqref{3.16} and the assumption $c_3=0$ we will be able to reduce the
exponent of $r$ in this inequality by $2$ every time, until we get
that $\tau$ must be bounded as $r\to0.$ As shown at step 1, then
$\tau$ is linearly dependent of $\sigma,$ which is a
contradiction. Thus $c_3\ne0.$ Then, \eqref{3.12} and \eqref{3.16}
give us $\tau'(r)=O(r^{-k})$ as $r\to0,$ or \eqref{3.13}.

From \eqref{3.13} and the equality
$\tau(r)=\tau(r_0)+\int_{r_0}^{r}\!\tau'(s)\,ds$ it follows that as
$r\to 0,$ $\tau(r)=O(\ln r)$ if $k=1,$ and $\tau(r)=O(r^{-k+1})$ if $k>
1.$ One can then get either \eqref{3.5a} or \eqref{3.5b} by conveniently
multiplying $\tau$ by a non-zero constant. This proves the lemma.
\end{proof}

The following lemma will study the properties of $j_m(r,l)$ for $l$ real.
\begin{lemma}\label{lm3.2}
Let $l=\lambda\in\mathbb R$. Then $j_m(r,\lambda)$ and $j'_m(r,\lambda)$ are real.
If $\lambda\le0,$ then $j_m(r,\lambda)\ge r^{|m|+1/2}$ and $j'_m(r,\lambda)>0.$
\end{lemma}
\begin{proof}
We will consider only the case $m\ge0,$ since $j_m(r,\lambda)$ is an
even function of $m.$ If $\lambda$ is real, then $q(r)-\lambda$ is
real. With the notation from the proof of lemma \ref{lm3.1},
$j_m(r,\lambda)=r^{m+1/2} \sigma(r,\lambda),$ with $\sigma(r,\lambda)$
given by the series \eqref{3.11}. By using this series, and
definition \eqref{3.7} of the operator $T$ we deduce that
$\sigma(r,\lambda)$ is real, and consequently, $j_m(r,\lambda)$ is
real.

If $\lambda\le 0,$ then $q(r)-\lambda\ge 0.$ In this case, the
operator $T$ will map non-negative functions into non-negative
functions. We conclude from the construction \eqref{3.11} of
$\sigma$ that $\sigma(r,\lambda)\ge 1,$ and so, $j_m(r,\lambda)\ge
r^{m+1/2}.$

To show $j'_m(r,\lambda)>0$ it suffices to prove
$\sigma'(r,\lambda)\ge0,$ which is an immediate consequence of
\eqref{3.9}.
\end{proof}

The next lemma will study the properties of $j_m(r,l)$ for large
absolute value of $m$ and $l=\lambda$ a real number.
\begin{lemma}\label{lm3.3}
Let $r>0$ be a fixed number and $\Lambda \subset [0,\infty)$ be a bounded set.
If $|m|$ is sufficiently large, then for all $\lambda\in \Lambda$
$$
  j_m(r,\lambda) > 0, \ \ j'_m(r,\lambda) > 0.
$$
\end{lemma}
\begin{proof}
We can assume $m\ge 0.$  Denote $k=2m+1\ge 1.$  We have that
$j_m(r,\lambda)=r^{m+1/2} \sigma(r,\lambda),$ with $\sigma(r,\lambda)$ defined by
the series \eqref{3.11}. Here we will prefer to use the notation $\sigma_m(r,\lambda)$
to emphasize its dependence on $m.$ By using this series and
estimate \eqref{3.8} it follows that
$$
  |1-\sigma_m(r,\lambda)|\le \frac{C}{k+1}, \ \ \lambda\in \Lambda,
$$
with $C=C(\Lambda,q_\infty)$ a constant, where $q_\infty=\sup_{r\in[0,\infty)}q(r).$

From \eqref{3.9} we deduce
$$
  |\sigma_m'(r,\lambda)|\le \frac{D}{k+1}, \ \ \lambda\in \Lambda,
$$
with $D=D(\Lambda,q_\infty)$ another constant.

These observations imply that for $m$ large enough and $\lambda\in
\Lambda$, $\sigma_m(r,\lambda)$ is close to 1, while $\sigma_m'(r,\lambda)$ is
close to zero.  Then, since
$$
  j_m(r,\lambda)=r^{m+1/2}\sigma_m(r,\lambda),
$$
and
$$
  j'_m(r,\lambda)=r^{m-1/2} \{(m+1/2)\, \sigma_m(r,\lambda)+r\,\sigma_m'(r,\lambda)\},
$$
it follows quickly that for $m$ large enough both of these quantities are strictly positive.
\end{proof}

\begin{lemma}\label{lm3.4}
For $l\in\{-\pi/2 < \arg(l-d^2) < 3\pi/2\},$ there exist
two solutions $w_m(r, l)$ and $x_m(r, l),$ $r\ge R, \, m\in\mathbb Z,$ of {\rm \eqref{2.5}}
which are analytic in $l.$ If $l\in \mathbb C_+=\{l: \mathfrak{Im}\,
l > 0\},$ then as $r\to\infty$, $w_m(r,l)$ together with its
derivative will decay exponentially, while $x_m(r,l)$ and its
derivative will increase exponentially.
\end{lemma}
\begin{proof}
Recall that $q(r)=d^2$ for $r\ge R,$ therefore on
$[R,\infty)$  \eqref{2.5} becomes
$$
  w''+\bigg\{l - d^2-\frac{m^2-1/4}{r^2}\bigg\}w=0.
$$
Set
\begin{equation}\label{3.17}
  w_m(r, l)=\sqrt{r}H^{(1)}_m(\sqrt{l-d^2}\,r), \ \
  x_m(r,l)=\sqrt{r}H^{(2)}_m(\sqrt{l-d^2}\,r),
\end{equation}
where $H^{(1)}_m(\zeta)$ and $H^{(2)}_m(\zeta)$ are the Hankel
functions of $m\mbox{-th}$ order.  They are analytic functions with
the domain $\{-\pi < \arg \zeta < \pi\}.$ Consider the function
$\sqrt{z}$ defined on $\{-\pi/2 < \arg z < 3\pi/2\}$ with values in
$\{-\pi/4 < \arg \zeta < 3\pi/4\}.$ We obtain that $w_m(r, l)$
and $x_m(r, l)$ are analytic functions of $l\in\{-\pi/2 < \arg(l-d^2)
< 3\pi/2\},$ which is the whole complex plane except those $l$ for
which $l-d^2$ has a zero real part and a non-positive imaginary part.

According to formulas (9.2.3) and (9.2.4) from,\cite{AS} we have the
asymptotic expansions
$$
  H^{(1)}_m(\zeta) \sim \sqrt{2/(\pi\zeta)}\, e^{i(\zeta-m\pi/2-\pi/4)}
$$
and
$$
  H^{(2)}_m(\zeta) \sim \sqrt{2/(\pi\zeta)}\, e^{-i(\zeta-m\pi/2-\pi/4)}
$$
as $|\zeta|\to\infty$ and $|\arg\zeta|<\pi.$
Then,
$$
  w_m(r,l)=O\big(e^{i\sqrt{l-d^2}\,r}\big),\ \
  x_m(r,l)=O\big(e^{-i\sqrt{l-d^2}\,r}\big) \mbox{ as } r\to\infty.
$$
If $l\in\mathbb C_+,$ then $\mathfrak{Im}\, \sqrt{l-d^2} > 0$ and
thus, $i\sqrt{l-d^2}$ has a strictly negative real part.
Therefore, as $r\to\infty,$ $w_m(r,l)$ will decay exponentially,
while $x_m(r,l)$ will increase exponentially. By formally
differentiating the above equalities (for a rigorous justification one
needs to use equalities (9.2.13) and (9.2.14) from \cite{AS}) we
deduce that $w_m'(r,l)$ will decay exponentially, while
$x_m'(r,l)$ will increase exponentially as $r\to\infty.$
\end{proof}

\section{Classification of the Solutions}\label{sec4}
The main purpose of this paper is to prove that under certain conditions,
the solution to the Helmholtz equation \eqref{2.1} is a superposition of
functions of the form \eqref{2.2}. For each of the functions in the
superposition, $v(r)$ will satisfy  \eqref{2.4} with the
notations of \eqref{2.3}, and $l$ will be a real variable, which we will
denote by $\lambda$ (thus we will reserve the notation $l$ for the
complex variable, and the notation $\lambda$ for its restriction to
the real axis). In addition the following properties will hold:
\begin{subequations}
  \begin{equation}\label{4.1a}
    v(r) \mbox{ is bounded as } r\to 0,
  \end{equation}
  and
  \begin{equation}\label{4.1b}
    r\to \sqrt{r}v(r)
\begin{cases}
      \mbox{ is in } L^2(R,\infty),       &\mbox{ if } \lambda \le d^2,\\
      \mbox{ is bounded as } r\to\infty,  &\mbox{ if } \lambda > d^2.
    \end{cases}
  \end{equation}
\end{subequations}
In this section we will study and classify the functions $v(r)$ with the
properties \eqref{4.1a} and \eqref{4.1b}.

A solution $v(r)$ of \eqref{2.4} has the form $v(r)=w(r)/\sqrt{r},$
with $w(r)$ satisfying \eqref{2.5}. Lemma \ref{lm3.1} shows how the solutions of
\eqref{2.5} look like. It is clear that in order that $v(r)$  satisfy \eqref{4.1a},
we need $w(r)=j_m(r,\lambda),$ or
\begin{equation}\label{4.2}
  v(r)=\frac{j_m(r,\lambda)}{\sqrt{r}}.
\end{equation}
Condition  \eqref{4.1b} is then satisfied if and only if
\begin{equation}\label{4.3}
  r\to j_m(r,\lambda)
\begin{cases}
    \mbox{ is in } L^2(R,\infty),       &\mbox{ if } \lambda \le d^2,\\
    \mbox{ is bounded as } r\to\infty,  &\mbox{ if } \lambda > d^2.
  \end{cases}
\end{equation}
Next we will investigate for which $\lambda$ condition \eqref{4.3}
holds.

We will need to consider four cases: $\lambda\le 0$, $0 < \lambda <
d^2,$ $\lambda=d^2$ and $\lambda>d^2.$ In each of these intervals
$j_m(r,\lambda)$ and consequently $v(r),$ will have a different
behavior.

{\bf Case 1.} If $\lambda\le 0,$ then according to lemma \ref{lm3.2},
$j_m(r,\lambda)\ge r^{|m|+1/2},$ therefore $j_m(r,\lambda)$ will be not square integrable
on $(R,\infty).$

{\bf Case 2.} Assume  $0 < \lambda < d^2.$ For $r\ge R$ the
function $q(r)$ defined in \eqref{2.3} is constant and equal to
$d^{2}.$ Then \eqref{2.5} becomes
\begin{equation}\label{4.4}
   w''+\bigg\{\lambda - d^2 - \frac{m^2-1/4}{r^2}\bigg\}w=0, \ \ r \in [R, \infty).
\end{equation}
The solutions to this equation are
\begin{equation}\label{4.5}
  k_m(r,\lambda)=\sqrt{r}K_m(\sqrt{d^2-\lambda}\,r), \ \ \lambda< d^2, r \ge R,
\end{equation}
and
$$
  \sqrt{r}I_m(\sqrt{d^2-\lambda}\,r), \ \ \lambda< d^2, r \ge R,
$$
where $K_m$ and $I_m$ are the modified Bessel functions. Formulas
(9.7.1) and (9.7.2) from \cite{AS} give us the expansions
\begin{equation*}
     \sqrt{s}K_m(s) \sim \sqrt{\pi/2}\,e^{-s} \mbox{ as } s\to\infty,
\end{equation*}
and
$$
  \sqrt{s}I_m(s) \sim \sqrt{\pi/2}\,e^{s} \mbox{ as } s\to\infty.
$$

Thus, we have one solution decaying exponentially while another increasing exponentially.
In order that $j_m(r,\lambda)$ be in $L^2(R,\infty)$ we  need
$$
  j_m(r,\lambda)=C k_m(r,\lambda) \mbox{ for } r\ge R,
$$
for some constant $C.$ Since both $j_m(r,\lambda)$ and $C
k_m(r,\lambda)$ satisfy the same second order differential equation,
namely \eqref{4.4}, to ask for the equality of these functions is the
same as to ask for them to satisfy the same boundary conditions at
$r=R,$ which translates into
$$
  j_m(R,\lambda)-C k_m(R,\lambda)=0,
$$
and
$$
  j'_m(R,\lambda)-C k'_m(R,\lambda)=0.
$$
Thus, $C$ must satisfy simultaneously two different conditions. This is possible if
and only if
\begin{equation}\label{4.6}
 \frac{j_m'(R, l)}{j_m(R, l)}=\frac{k_m'(R, l)}{k_m(R, l)}, \ \ \lambda < d^2.
\end{equation}
Then we will have
\begin{equation}\label{4.7}
  j_m(r,\lambda)=\frac{j_m(R,\lambda)}{k_m(R,\lambda)}k_m(r,\lambda),\ \ r\ge R.
\end{equation}

It will be shown later that for each $m,$ the set of $\lambda$ such
that \eqref{4.6} holds is finite. Notice that in this case
$j_m(r,\lambda)$ will decay exponentially as $r\to\infty.$ Its
derivative has the same property, this follows from the asymptotic
expansion
\begin{equation}\label{4.8}
   [\sqrt{s}K_m(s)]' \sim - \sqrt{\pi/2} \,e^{-s} \mbox{ as } s\to\infty
\end{equation}
(according to the formulas (9.7.2) and (9.7.4) from.\cite{AS})
In conclusion, the only $\lambda\in(0, d^2)$ for which \eqref{4.3} holds,
are those satisfying \eqref{4.6}.

{\bf Case 3}. Let now $\lambda=d^2.$
Two linear independent solutions of \eqref{4.4} are in this case
$$
  r \to r^{1/2-|m|} , \ \ m\in\mathbb Z,
$$
and
$$
 r\to
\begin{cases}
         \sqrt{r} \ln r & \mbox{ for } m=0,  \\
         r^{1/2+|m|}    & \mbox{ for } m\ne 0.
      \end{cases}
$$
The second solution is not bounded for any $m\in\mathbb Z.$ By
matching the boundary conditions at $r=R$ as in the previous case, it
is easy to show that $j_m(r,\lambda)$ will be proportional to the
first of these two solutions if and only if
\begin{equation}\label{4.9}
  \frac{j_m'(R, l)}{j_m(R, l)}=-\frac{|m|-1/2}{R}, \ \ \lambda=d^2,
\end{equation}
and then,
\begin{equation}\label{4.10}
  j_m(r,\lambda)=\frac{j_m(R,\lambda)}{R^{1/2-|m|}}r^{1/2-|m|},\ \ r\ge R.
\end{equation}
The function $j_m(r,\lambda)$ will be in $L^2(R,\infty)$ if and only if $|m|\ge 2.$

A function of the form \eqref{2.2}, with $v(r)$ given by
\eqref{4.2}, for which $0<\lambda\le d^2$ and either
\eqref{4.6}, or \eqref{4.9} (with $|m|\ge 2$) holds, is
called a \emph{guided mode}. Note that a guided mode decays in $r$
either exponentially, or as $r^{-|m|}$ $(|m|\ge 2).$

{\bf Case 4.} Let $\lambda>d^2.$ Two solutions of \eqref{4.4}
are then
\begin{equation}\label{4.11}
  a_m(r,\lambda)=\sqrt{r}J_m\big(\sqrt{\lambda-d^2}\,r\big), \ \
  b_m(r,\lambda)=\sqrt{r}Y_m\big(\sqrt{\lambda-d^2}\,r\big),  \ \  \lambda > d^2,
\end{equation}
where $J_m$ and $Y_m$ are the Bessel functions of the first and second kind.
Note the formulas
$$
  \sqrt{s}\,J_m(s)=\sqrt{\pi/2}\cos(s-m\pi/2-\pi/4)+O(s^{-1/2}),
$$
$$
  \sqrt{s}\,Y_m(s)=\sqrt{\pi/2}\sin(s-m\pi/2-\pi/4)+O(s^{-1/2}),
$$
as $s\to\infty$ (they are a particular case of formulas (9.2.1) and
(9.2.2) from.\cite{AS}) We infer that $a_m(r,\lambda)$ and
$b_m(r,\lambda)$ will be bounded as $r\to \infty.$ By formally
differentiating the above formulas it follows that $a'_m(r,\lambda)$
and $b'_m(r,\lambda)$ will also be bounded as $r\to\infty.$

The functions $a_m(r,\lambda)$ and $b_m(r,\lambda)$  are linearly independent,
since $J_m$ and $Y_m$ are  linearly independent. Then, for $r\ge R,$
$j_m(r,\lambda)$ will be a linear combination of them.  To find the coefficients
of the linear combination, set
$$
  j_m(r,\lambda)=c_m(\lambda)a_m(r,\lambda)+d_m(\lambda)b_m(r,\lambda), \ \ r \ge R.
$$
By matching the boundary conditions at $r=R$, we obtain a linear
system which enables us to solve for $c_m$ and $d_m.$ Apply the
equality
$$
  Y'_m(z)J_m(z)-J'_m(z)Y_m(z)=2/(\pi z),
$$
((9.1.16) from \cite{AS}) to find a value for the determinant of the this linear system.
Obtain
\begin{equation}\label{4.12}
  b_m'(R,\lambda)a_m(R,\lambda)-a_m'(R,\lambda)b_m(R,\lambda)=2/\pi,
\end{equation}
and therefore,
\begin{subequations}
\begin{equation}\label{4.13a}
  c_m(\lambda)=\frac{\pi}{2}\{b_m'(R,\lambda)j_m(R,\lambda)-j_m'(R,\lambda)b_m(R,\lambda)\},
\end{equation}
\begin{equation}\label{4.13b}
  d_m(\lambda)=-\frac{\pi}{2}\{a_m'(R,\lambda)j_m(R,\lambda)-j_m'(R,\lambda)a_m(R,\lambda)\}.
\end{equation}
\end{subequations}

It is easy to see that $j_m(r,\lambda)$ and its derivative will be
bounded as $r\to \infty,$ and so, condition \eqref{4.3} will be
satisfied for all $\lambda>d^2.$

Let us look at the expression of \eqref{2.2} for $\lambda>d^2$
and with $v(r)$ given by \eqref{4.2}. Notice that if $d^2 < \lambda <
k^2 n_0^2,$ then \eqref{2.2} will be oscillatory in $z$
(recall that $k^2\beta^2=k^2n_0^2-\lambda$). In this case we will say
that \eqref{2.2} is a \emph{radiation mode}.  On the other hand,
if $\lambda > k^2 n_0^2,$ then $\beta$ becomes imaginary. Depending
on the sign of $\mathfrak{Im} \beta$ we will have exponential decay in
one of the directions $z\to-\infty,$ $z\to\infty,$ and exponential
growth in the other one. For $\lambda > k^2 n_0^2,$ \eqref{2.2}
will be called an \emph{evanescent mode}.

\section{The Theory of Eigenvalue Problems}\label{sec5}
Consider the eigenvalue problem
\begin{equation}\label{5.1}
  w''+\{l-Q(r)\}w=0, \ \ r\in (0,\infty),
\end{equation}
where $l\in\mathbb C.$ Assume that $Q$ is integrable over any compact
subset of $(0, \infty)$ (in \cite{CL} and \cite{Ti} the theory is
developed only for continuous functions $Q,$ but it is mentioned in a
footnote on page 224 of \cite{CL} that it suffices for $Q$ to be as we
assume above). Let $0<R<\infty$ be an arbitrary but fixed number.  Let
$\varphi(r,l)$ and $\theta(r,l)$ be the solutions of \eqref{5.1} with
the boundary conditions
\begin{equation}\label{5.2}
\begin{cases}
  \varphi(R, l)  = 0, \ \ \varphi'(R, l)=-1,\\
  \theta(R, l)   = 1, \ \ \theta'(R, l)=0.
\end{cases}
\end{equation}
Since \eqref{5.1} has an analytic dependence on the parameter $l,$ the
solutions $\theta(r,l)$ and $\varphi(r,l)$ will be analytic functions
of $l$ for $r$ fixed.

Any solution to \eqref{5.1} linearly independent of $\varphi$ can be
represented, up to a constant multiple, in the form
\begin{equation}\label{5.3}
\psi=\theta+M\varphi,
\end{equation}
with $M\in\mathbb C.$ Let $\tau\in\mathbb R$ and $0<t<\infty.$
Look for a solution of \eqref{5.1} of the form \eqref{5.3} to satisfy
the boundary condition
$$
  \cos\tau\psi(t,l)+\sin\tau\psi'(t,l)=0.
$$
It is a direct calculation to check that we need
\begin{equation}\label{5.4}
M=M(l)=-\frac{\theta(t,l)\cos\tau+\theta'(t,l)\sin\tau}
                          {\varphi(t,l)\cos\tau+\varphi'(t,l)\sin\tau}.
\end{equation}

Let $\mathbb C_+$ be the open upper complex-half-plane, $\mathbb
C_+=\{l: \mathfrak{Im}\, l > 0\}.$ As shown in chapter 2 of \cite{Ti}
the following results hold: as $t\to0,$ $M(l)$ converges uniformly on
compact subsets of $\mathbb C_+$ to a function $M_0(l)$ analytic on
$\mathbb C_+.$ Moreover, the function
\begin{equation}\label{5.5}
  \psi_0(r, l)= \theta(r, l)+M_0(l)\,\varphi(r, l), \ \ l\in \mathbb C_+,
\end{equation}
is in $L^2(0, R),$ and one has
\begin{equation}\label{5.6}
    \int\limits_{0}^{R}     \!|\psi_0(r,l)|^2 \,dr =  \frac{\mathfrak{Im}\, M_0(l)}{\mathfrak{Im}\, l}.
\end{equation}
As $t\to\infty,$ $M(l)$ converges uniformly on compact subsets of
$\mathbb C_+$ to an analytic function $M_{\infty}(l)$ on $\mathbb C_+,$
and if
\begin{equation}\label{5.7}
  \psi_\infty(r, l)= \theta(r, l)+M_{\infty}(l)\,\varphi(r, l), \ \
  l\in \mathbb C_+,
\end{equation}
then  $\psi_\infty(r, l)\in L^2(R, \infty)$ and
\begin{equation}\label{5.8}
  \int\limits_{R}^{\infty}\!|\psi_\infty(r,l)|^2 \,dr = -\frac{\mathfrak{Im}\, M_{\infty}(l)}{\mathfrak{Im}\, l}.
\end{equation}

Note that the obtained $M_0, \, M_{\infty},\, \psi_0,\, \psi_\infty$ depend
on the parameter $\tau\in\mathbb R.$ Thus, possibly these quantities,
and therefore the representation given below, in Theorem \ref{th5.1},
will not be unique. This might be true in general, but in our concrete
case, given by  \eqref{2.5}, these quantities will turn out to
be unique, as we will see from lemma \ref{lm6.1}. So then the transform we are
looking for (which is calculated in Theorems \ref{th7.1} and
\ref{th7.2}) will be unique.

In section 9.5 of \cite{CL} and chapter 3 of \cite{Ti} it is proved
that for any $\lambda\in \mathbb R$ the following limits exist
\begin{subequations}
\begin{equation}\label{5.9a}
  \xi(\lambda)=\lim\limits_{\delta\to0^+}\int\limits_0^\lambda\!
  -\mathfrak{Im}\frac{1}{M_0(s+i\delta)-M_{\infty}(s+i\delta)}\,d s,
\end{equation}
\begin{equation}\label{5.9b}
  \eta(\lambda)=\lim\limits_{\delta\to0^+}\int\limits_0^\lambda\!
  -\mathfrak{Im}\frac{M_0(s+i\delta)}{M_0(s+i\delta)-M_{\infty}(s+i\delta)}\,d s,
\end{equation}
\begin{equation}\label{5.9c}
  \zeta(\lambda)=\lim\limits_{\delta\to0^+}\int\limits_0^\lambda\!
  -\mathfrak{Im}\frac{M_0(s+i\delta)M_{\infty}(s+i\delta)}{M_0(s+i\delta)-M_{\infty}(s+i\delta)}\,d s.
\end{equation}
\end{subequations}
It is shown there that the functions $\xi$ and $\zeta$ are
non-decreasing, and that $\eta$ is with bounded variation.  In
addition, for any $\lambda_0 < \lambda_1$ real numbers,
\begin{equation}\label{5.10}
  \{\eta(\lambda_1)-\eta(\lambda_0)\}^2 \le
  \{\xi(\lambda_1)-\xi(\lambda_0)\}\{\zeta(\lambda_1)-\zeta(\lambda_0)\},
\end{equation}
as stated on page 252 of \cite{CL} (with a different notation).  Then,
equalities (3.1.8), (3.1.9), (3.1.10) from \cite{Ti} give us an
expansion formula for a function $g\in L^{2}(0,\infty)$ in terms of
$\theta(r,l),\,\varphi(r,l)$ and the functions
$\xi,\,\eta$ and $\zeta.$ The same result is proved in \cite{CL} at
Theorem 5.2. In this reference the representation result is stated
more rigorously, so we will prefer it over.\cite{Ti} To state
the result we need some notations.

Denote $\rho=(\xi,\,\eta,\,\zeta).$ For any vector $\Gamma=(\Gamma_1,\Gamma_2),$
where $\Gamma_1,\,\Gamma_2\,:\mathbb R\to \mathbb C,$ let
\begin{equation}\label{5.11}
  ||\Gamma||^2=\int\limits_{-\infty}^{\infty}\! |\Gamma_1(\lambda)|^2\, d\xi
              +2{\mathfrak{Re}}\{\Gamma_1(\lambda)\,\bar \Gamma_2(\lambda)\}\,d\eta
              + |\Gamma_2(\lambda)|^2\,d\zeta.
\end{equation}
The fact that $\xi$ and $\eta$ are non-decreasing, together with
\eqref{5.10},  gives us that $||\Gamma||^2\ge 0.$ It is  easy to
check that $||\cdot||$ is a semi-norm.  Denote by $L^{2}(\rho)$ the
space of all $\Gamma=(\Gamma_1, \Gamma_2)$ such that $||\Gamma||<\infty.$ This is then the
statement of Theorem 5.2 from.\cite{CL}

\begin{theorem}\label{th5.1}
If $g\in L^2(0,\infty),$ the vector $\Gamma=(\Gamma_1,\Gamma_2)$, where
$$
  \Gamma_1(\lambda)=\int\limits_{0}^{\infty}\!\theta(r,\lambda)g(r)\,dr, \ \
  \Gamma_2(\lambda)=\int\limits_{0}^{\infty}\!\varphi(r,\lambda)g(r)\,dr,
$$
converges in $L^{2}(\rho),$ that is, there exists $\Gamma\in
L^{2}(\rho)$ such that
$$
  ||\Gamma-\Gamma^{cd}||\to 0 \mbox{ as } c\to 0, \, d\to \infty,
$$
where for $0 < c< d < \infty$
\begin{equation}\label{5.12}
  \Gamma^{cd}_1(\lambda)=\int\limits_{c}^{d}\!\theta(r,\lambda)g(r)\,dr, \ \
  \Gamma^{cd}_2(\lambda)=\int\limits_{c}^{d}\!\varphi(r,\lambda)g(r)\,dr.
\end{equation}
The expansion
\begin{multline*}
 g(r)= \frac{1}{\pi} \int\limits_{-\infty}^{\infty}\!\big\{
             \theta(r,\lambda)\Gamma_1(\lambda)\,d\xi(\lambda)+\theta(r,\lambda)\Gamma_2(\lambda)\,d\eta(\lambda)\\
            +\varphi(r,\lambda)\Gamma_1(\lambda)\,d\eta(\lambda)+\varphi(r,\lambda)\Gamma_2(\lambda)\,d\zeta(\lambda)\big\}
\end{multline*}
holds, with the latter integral convergent in $L^{2}(0,\infty),$ that
is, $g^{\sigma\tau}\to g$ in $L^{2}(0,\infty)$ as $\sigma\to-\infty,\,
\tau\to\infty,$ where for $-\infty < \tau < \sigma <\infty$
\begin{multline}\label{5.13}
  g^{\sigma\tau}(r)= \frac{1}{\pi} \int\limits_{\sigma}^{\tau}\!\big\{
     \theta(r,\lambda)\Gamma_1(\lambda)\,d\xi(\lambda)+\theta(r,\lambda)\Gamma_2(\lambda)\,d\eta(\lambda) \\
    +\varphi(r,\lambda)\Gamma_1(\lambda)\,d\eta(\lambda)+\varphi(r,\lambda)\Gamma_2(\lambda)\,d\zeta(\lambda)\big\}.
\end{multline}
We have the Parseval identity
\begin{equation}\label{5.14}
  \int\limits_{0}^{\infty}\!|g(r)|^2\,dr=\frac{1}{\pi}||\Gamma||^2.
\end{equation}
\end{theorem}

\section{Computing the Measures}\label{sec6}
In the sections to follow we will apply the results of section 5 to
our particular eigenvalue equation, given by \eqref{2.5}. Thus, for the function
$Q(r)$ in  \eqref{5.1} we will have the expression
$$
  Q(r)=q(r) - \frac{m^2-1/4}{r^2}, \ \ r\in (0, \infty).
$$
In this section we will calculate the measures given in \eqref{5.9a},
\eqref{5.9b} and \eqref{5.9c}.

\begin{lemma}\label{lm6.1}
Let $M_0^m(l),\, M_{\infty}^m(l), \psi^m_0(r,l),\, \psi^m_\infty(r,l)$
$(l\in\mathbb C_+=\{l: \mathfrak{Im}\, l >0\})$ be the quantities
defined in section {\rm 5} for equation {\rm \eqref{2.5}}
(we use the superscript $m$ to emphasize their dependence on $m\in\mathbb Z$).
Let $j_m(r,l)$ $w_m(r,l)$ be the solutions of {\rm \eqref{2.5}} defined respectively in
lemma \ref{lm3.1} and  by {\rm \eqref{3.17}}. Then,
\begin{equation}\label{6.1}
  M_0^m(l)=-\frac{j_m'(R, l)}{j_m(R, l)}, \ \
  M_{\infty}^m(l)=-\frac{w_m'(R, l)}{w_m(R, l)},
\end{equation}
and
\begin{equation}\label{6.2}
  \psi^m_0(r, l)=\frac{j_m(r, l)}{j_m(R, l)}, \ \
  \psi^m_\infty(r, l)=\frac{w_m(r, l)}{w_m(R, l)}.
\end{equation}
\end{lemma}
\begin{proof}
It is easy to note that $j_m(r, l)$ and $y_m(r, l)$ defined in lemma \ref{lm3.1}
are linearly independent solutions of \eqref{2.5}, thus any
other solution will be a linear combination of these two. Then $\theta(r,l)$
and $\phi(r,l)$ defined in section 5 can be represented as
\begin{equation}\label{6.3}
\begin{cases}
  \theta(r,l)=\alpha(l)j_m(r, l)+\beta(l)y_m(r, l),\\
  \varphi(r,l)=\gamma(l)j_m(r, l)+\delta(l)y_m(r, l),
\end{cases}
\end{equation}
for some coefficients $\alpha,\,\beta,\,\gamma,\delta.$ Let
\begin{equation}\label{6.4}
  \Delta(l)=j_m'(R,l)y_m(R,l)-j_m(R,l)y_m'(R,l).
\end{equation}
Being the Wronskian of two linearly independent solutions $\Delta(l)$ is
non-zero. Using the boundary conditions \eqref{5.2} it is
easy to find that
\begin{equation}\label{6.5}
  \alpha=-\frac{y_m'(R, l)}{\Delta(l)},\, \beta=\frac{j_m'(R,
    l)}{\Delta(l)},\ \gamma=-\frac{y_m(R, l)}{\Delta(l)},\,
  \delta=\frac{j_m(R, l)}{\Delta(l)}.
\end{equation}

From \eqref{6.3} and \eqref{5.4} it follows that
$$
  M(l)=-\frac{\{\alpha j_m(t, l)+\beta y_m(t, l)\}\cos\tau
  +\{\alpha j_m'(t, l)+\beta y_m'(t, l)\}\sin\tau}
  {\{\gamma j_m(t, l)+\delta y_m(t, l)\}\cos\tau
  +\{\gamma j_m'(t, l)+\delta y_m'(t, l)\}\sin\tau}.
$$
For $t\to0$ the formulas \eqref{3.1}, \eqref{3.2a} and \eqref{3.2b} tell us
that the term $y_m(t,l)$ will dominate $j_m(t,l)$ and
$j_m'(t,l)$, while $y_m'(t,l)$ will dominate $y_m(t,l).$
Then,
$$
  M_0^m(l)=\lim\limits_{t\to 0} M(l)=-\frac{\beta}{\delta}.
$$
By applying \eqref{6.5} we obtain the first equality in \eqref{6.1}. To get the first
equality in \eqref{6.2} use definition \eqref{5.5} of $\psi^m_0(r,l)$ and
equalities \eqref{6.3} to \eqref{6.5}.

The quantities $M_{\infty}^m(l)$ and $\psi^m_\infty(r,l)$ are calculated
in exactly the same way.  One needs to express $\varphi(r, l)$ and
$\theta(r,l)$ in terms of $w_m(r,l)$ and $x_m(r,l),$ then put
$t\to\infty$ in \eqref{5.4} and use the properties of $w_m(r,l)$
$x_m(r,l),$ and their derivatives shown in lemma \ref{lm3.4}.
The lemma is proved.
\end{proof}

The next lemma will study the properties of the functions $M_0^m(l),$
$M_{\infty}^m(l),$ and $M_0^m(l)-M_{\infty}^m(l).$ Here we will set some notation and write
some formulas which will be used in the lemma. By $\lambda$ we
will denote a real variable.  Let $w_m(r,\lambda),$ $k_m(r,\lambda),$ $a_m(r,\lambda),$
and $b_m(r,\lambda)$ be the functions defined by \eqref{3.17}, \eqref{4.5} and
\eqref{4.11}.  The following equalities hold,
$$
  H^{(1)}_m(iz)=\frac{2}{i\pi}e^{-im/2}K_m(z), \ \ -\pi < \arg z  \le\pi/2
$$
and
\begin{equation}\label{6.6}
  H^{(1)}_m(z)=J_{m}(z)+iY_{m}(z), \ \ -\pi < \arg z < \pi,
\end{equation}
(formulas  (9.6.4) and (9.1.3) from.\cite{AS}) We deduce
\begin{subequations}
  \begin{equation}\label{6.7a}
    w_m(r,\lambda)=\frac{2}{i\pi}e^{-im/2}k_m(r,\lambda), \ \ \lambda < d^2,
  \end{equation}
  and
  \begin{equation}\label{6.7b}
    w_m(r,\lambda)=a_m(r,\lambda)+i b_m(r,\lambda)  \ \ \lambda > d^2.
  \end{equation}
\end{subequations}
\begin{lemma}\label{lm6.2}
$M_0^m(l)$ is meromorphic across all the complex plane, while
$M_{\infty}^m(l)$ $(l\in\mathbb C_+)$ extends continuously to the real
axis.  $M_0^m(\lambda)-M_{\infty}^m(\lambda)$ is real or infinite for
$\lambda < d^2,$ and it has a finite number of zeros on the real axis,
all in the interval $(0, d^2].$
\end{lemma}
\begin{proof} As stated in lemma \ref{lm3.1}, $j_m(R,l)$ and
$j_m'(R,l)$ will be analytic functions of $l\in\mathbb C.$ Then
$M_0^m(l),$ being obtained as their ratio, will be a meromorphic
function.

The function $w_m(r,l)$ was defined in lemma \ref{lm3.4}. It was
proved there that $w_m(r,l)$ is defined and analytic for all
$l\in\{-\pi/2 < \arg(l-d^2) < 3\pi/2\}.$ In particular
$w_m(r,\lambda)$ is defined for all real $\lambda\ne d^{2}.$ To
prove that $M_{\infty}^m(l)$ extends continuously to the real
axis, it suffices to show that its denominator, $w_m(r,\lambda),$
is not zero for $\lambda\in \mathbb R\backslash\{d^2\}$ and that
$\lim_{l\to d^2} M_{\infty}^m(l)$ exists.

If $\lambda<d^2,$ then $w_m(r,\lambda)\ne0$ because of \eqref{4.5} and
\eqref{6.7a}, since the Bessel function $K_m(s)$ takes real strictly
positive values for $s>0$ (as it follows from section (9.6.1) of.\cite{AS})
If $\lambda>d^2,$ then $w_m(r,\lambda)\ne0$ because
of \eqref{4.11} and \eqref{6.7b}, since for $s>0$ the Bessel functions
$J_m(s)$ and $Y_m(s)$ take real values and cannot be zero at the same
time.

Show that $\lim_{l\to d^2} M_{\infty}^m(l)$ exists. If $m\ge 0,$ then according to
(9.1.8) and (9.1.9) from, \cite{AS} we have that for $z\in\mathbb C$,
$z\to0$
$$
  H^{(1)}_m(z)\sim -(1/\pi)(m-1)!(z/2)^{-m}, \mbox{ for }m\ge 1,
$$
and
$$
  H^{(1)}_m(z)\sim (-2i/\pi)\ln{z},  \mbox{ for }m=0.
$$
We can extend these to $m <0,$ by using \eqref{6.6} together with
$$
  J_{-m}(z)=(-1)^{m}J_m(z), Y_{-m}(z)=(-1)^{m}Y_m(z), \ \ m\in \mathbb Z
$$
(formula (9.1.5) from.\cite{AS}) One can obtain the behavior of the
derivative of $H^{(1)}_m(z)$ as $z\to0$ by formally differentiating
the above. Set $z=R\sqrt{l-d^2}.$ Deduce that for all $m\in\mathbb Z$
\begin{equation}\label{6.8}
  \lim\limits_{l\to d^2} M_{\infty}^m(l)=\frac{|m|-1/2}{R}.
\end{equation}

Note that $M_0^m(\lambda)$ is real or infinite for $\lambda< d^2,$ being the
quotient of $-j_m'(R,\lambda)$ and $j_m(R,\lambda)$ both of which are real
according to lemma \ref{lm3.1}. We have that $M_{\infty}^m(\lambda)$ is real for
$\lambda< d^2,$ that follows from \eqref{6.7a} and by using again the
properties of the function $K_m(s).$ Then, $M_0^m(\lambda)-M_{\infty}^m(\lambda)$ is real
or infinite for $\lambda< d^2.$

Let us show the last part of this lemma, the fact that
$M_0^m(\lambda)-M_{\infty}^m(\lambda)$ finite number of zeros on the real axis, all in
the interval $(0, d^2].$ Let first prove that this function can have
no zeros for $\lambda\le0.$ From \eqref{6.7a} we have
$$
  M_0^m(\lambda)-M_{\infty}^m(\lambda)
  =-\frac{j_m'(R,\lambda)}{j_m(R,\lambda)}+\frac{k_m'(R,\lambda)}{k_m(R,\lambda)}
  =-\frac{D_m(R,\lambda)}{j_m(R,\lambda)k_m(R,\lambda)},
$$
where
$$
  D_m(r,\lambda)=j_m'(r,\lambda)k_m(r,\lambda)-j_m(r,\lambda)k_m'(r,\lambda).
$$

Assume for some $\lambda,$ $M_0^m(\lambda)-M_{\infty}^m(\lambda)=0.$ Then $D_m(R,\lambda)=0$. We will
prove that  is false.  First note that $D_m(r,\lambda)=D_m(R,\lambda)$ for
$r\ge R.$ Indeed, $j_m(r,\lambda)$ and $k_m(r,\lambda)$ will satisfy \eqref{2.5}
(for $j_m(r,\lambda)$ this follows by its definition, for $k_m(r,\lambda)$ it
follows from the observation that according to \eqref{6.7a}, $k_m(r,\lambda)$
is a constant times $w_m(r,\lambda)$ which, as defined by \eqref{3.17}, is a
solution of \eqref{2.5} for $r\ge R$). Then it is easy to check that
$$
  D_m'(r, \lambda)=j_m''(r,\lambda)k_m(r,\lambda)-j_m(r,\lambda)k_m''(r,\lambda)=0,
$$
 so $D_m(r, \lambda)$ is constant in $r.$ Second, note that for $r$
large enough $D_m(r,\lambda) < 0.$ Indeed, for $r > 0,$
$j_m(r,\lambda) >0,\,j_m'(r,\lambda) >0$ by lemma \ref{lm3.2},
and $k_m(r,\lambda)>0$ by the properties of Bessel functions.
Also, for $r$ sufficiently large we have $k_m'(r,\lambda)<0,$ this
is a consequence of \eqref{4.8}. We infer that for $r$ large
enough $D_m(r,\lambda) < 0.$ These two observations imply
$D_m(R,\lambda) < 0,$ therefore $D_m(R,\lambda)$ is non-zero, and
$M_0^m(\lambda)-M_{\infty}^m(\lambda)\ne0.$

Now let $0<\lambda\le d^2.$ We showed that $M_0^m(l)$ is meromorphic
on the whole complex plane. The function $w_m(r,l),$ as defined by
\eqref{3.17}, is analytic on $\{-\pi/2 < \arg(l-d^2) < 3\pi/2\},$
in particular it is analytic on the set of $l$ such that ${\mathfrak{Re}}\,l
< d^2.$ Then $M_{\infty}^m(l)=-w_m'(R,l)/w_m(R,l)$ is
meromorphic on the same region, and so is $M_0^m(l)-M_{\infty}^m(l).$
Therefore, it can have only a discrete number of zeros on the interval
$(0, d^2).$

Assuming that the number of zeros is infinite, their only possible
accumulation point is $\lambda=d^2.$ So, there will exist a sequence
$\lambda_n < d^2, \, n \ge 1,$ with $\lim \lambda_n=d^2$ and
$M_0^m(\lambda_n)-M_{\infty}^m(\lambda_n)=0$ for all $n\ge 1.$ According to  formulas
(9.1.3), (9.1.10) and (9.1.11) from, \cite{AS} the Hankel function
$H^{(1)}(z)$ will have the representation
$$
  H^{(1)}_m(z)=f_1(z)+f_2(z)\ln{z},
$$
with $f_1(z)$ and $f_2(z)$ meromorphic functions of $z\in \mathbb
C.$ Then, if we recall  definition \eqref{3.17} of $w_m(r,l),$
one can calculate that we will have the representation
$$
  M_0^m(l)-M_{\infty}^m(l)=\frac{g_1(t)+g_2(t)\ln{t}}{g_3(t)+g_4(t)\ln{t}},
$$
with $g_1, g_2, g_3$ and $g_4$ meromorphic functions on the whole complex plane, and
$t=\sqrt{l-d^2}\in \mathbb C.$ We infer
$$
  g_1(t_n)+g_2(t_n)\ln{t_n}=0,
$$
for all $n\ge 1,$ where $t_n=\sqrt{\lambda_n-d^2}.$ Two cases are
possible.  If $g_2(t_n)$ is zero for infinitely many $n,$ this implies
$g_1(t_n)=0$ at the same points. Then the meromorphic functions
$g_1(t)$ and $g_2(t)$ are identically zero, and so is
$M_0^m(l)-M_{\infty}^m(l),$ obtaining a contradiction.  Otherwise, if
$g_2(t_n)$ is zero only for finitely many $n,$ we can write
$$
  \ln{t_n}=-\frac{g_1(t_n)}{g_2(t_n)}, \ \ n \ge n_0.
$$
On the right-hand side we have a meromorphic function. As $n\to\infty$ we have $t_n\to 0,$ thus, for some integer
$p$, there follows
$$  \ln|t_n| = O(|t_n|^p) .$$
That is clearly impossible. The contradiction shows that $M_0^m(\lambda)-M_{\infty}^m(\lambda)$ can have only
finitely many zeros on $[0, d^2).$ Another possible zero could be at $\lambda=d^2.$

Lastly, show that $M_0^m(\lambda)-M_{\infty}^m(\lambda)$ is never
zero on $(d^{2},\infty).$ According to \eqref{6.1} and \eqref{6.7b},
$$
  M_{\infty}^m(\lambda)=-\frac{a_m'(R,\lambda)+ib_m'(R,\lambda)}{a_m(R,\lambda)+ib_m(R,\lambda)},
$$
and so,
\begin{equation}\label{6.9}
\begin{split}
  M_0^m(\lambda)-M_{\infty}^m(\lambda) & =-\frac{j_m'}{j_m}+\frac{a_m'+ib_m'}{a_m+ib_m}\\
                                       & =-\frac{(j_m'a_m-a_m'j_m)+i(j_m'b_m-b_m'j_m)}{j_m(a_m+ib_m)}
\end{split}
\end{equation}
(the arguments $(R,\lambda)$ were omitted for simplicity).
If we assume that for some $\lambda> d^2,$ $M_0^m(\lambda)-M_{\infty}^m(\lambda)=0,$ this
implies
\begin{align*}
&  j_m'(R,\lambda)a_m(R,\lambda)-a_m'(R,\lambda)j_m(R,\lambda)=0, \\
&  j_m'(R,\lambda)b_m(R,\lambda)-b_m'(R,\lambda)j_m(R,\lambda)=0
\end{align*}
(since these quantities are real, as it follows from lemma \ref{lm3.2}
and \eqref{4.11}).  The numbers $j_m(R,\lambda)$ and $j_m'(R,\lambda)$
cannot be both zero, since $j_m(r,\lambda)$ is a non-zero solution of
the second order differential equation \eqref{2.5}. In this case the vectors
$(a_m(R,\lambda), a_m'(R,\lambda))$ and $(b_m(R,\lambda),
b_m'(R,\lambda))$ must be linearly dependent. That cannot be since the
functions $a_m(r,\lambda)$ and $b_m(r,\lambda)$ $(r\ge R)$ are also
solutions of \eqref{2.5}, and they are linearly independent.  Thus
$M_0^m(\lambda)-M_{\infty}^m(\lambda)\ne0.$ This finishes the lemma.
\end{proof}

\begin{theorem}\label{th6.3}
Let $\xi_m,$ $\eta_m$ and $\zeta_m$ be the functions defined by {\rm \eqref{5.9a},
\eqref{5.9b}} and {\rm \eqref{5.9c}}.  There exists a non-decreasing function
$\chi_m:\mathbb R\to\mathbb R$ such that the following measures are
equal
\begin{equation}\label{6.10}
\begin{split}
  d\xi_m(\lambda)   &=   j_m(R,\lambda)^2\,d\chi_m(\lambda), \\
  d\eta_m(\lambda)  &=  -j_m'(R,\lambda)j_m(R,\lambda)\,d\chi_m(\lambda), \\
  d\zeta_m(\lambda) &=   j_m'(R,\lambda)^2\,d\chi_m(\lambda).
\end{split}
\end{equation}
The function $\chi_m$ is identically zero for $\lambda\in (-\infty, 0],$ is
piecewise constant for $\lambda\in (0, d^2)$ where it has a finite number
of discontinuities, and is continuous for $\lambda\in (d^2, \infty).$
\end{theorem}

\begin{figure}[t]\label{fig3}
\begin{center}
\includegraphics[width=0.8\textwidth]{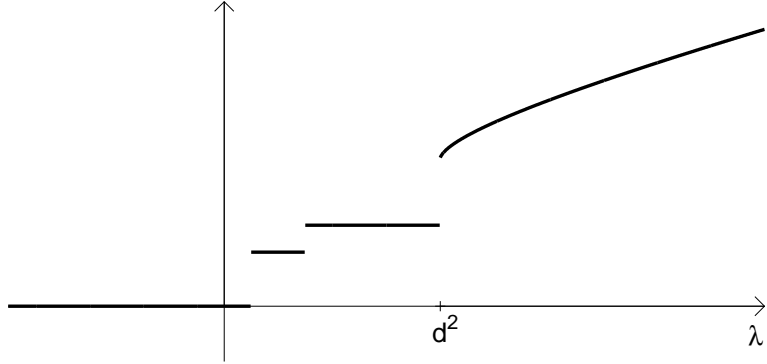}
\caption{The function $\chi_m(\lambda)$.}
\end{center}
\end{figure}

\begin{proof} Denote
$$
  M(l)=-\frac{1}{M_0^m(l)-M_{\infty}^m(l)}, \ \ l\in \mathbb C_+.
$$
According to \eqref{5.9a}, for any $\lambda_0 < \lambda_1$ real numbers
$$
  \xi_m(\lambda_1)-\xi_m(\lambda_0)
  =\lim\limits_{\delta\to0^+}\int\limits_{\lambda_0}^{\lambda_1}\! \mathfrak{Im}\, M(s+i\delta)\,d s.
$$
In particular, if $M(l)$ extends continuously to the interval
$[\lambda_0, \lambda_1],$ then by using Lebesgue's theorem of dominant
convergence it is easy to show that
\begin{equation}\label{6.11}
  \xi_m(\lambda_1)-\xi_m(\lambda_0)=\int\limits_{\lambda_0}^{\lambda_1}\!  \mathfrak{Im}\, M(s) \,d s.
\end{equation}
The same kind of reasoning clearly holds for $\eta_m$ and $\zeta_m.$

As it follows from lemma \ref{lm6.2}, if $\lambda \le 0$ then
$M_0^m(\lambda)-M_{\infty}^m(\lambda)$ is real or infinite, and
non-zero. By applying \eqref{6.11} we find that for any $\lambda_1 <
\lambda_2 < 0,$ $\xi_m(\lambda_1)-\xi_m(\lambda_0)=0.$ From \eqref{5.9a}
we have that $\xi_m(0)=0,$ and thus $\xi_m(\lambda)=0$ for all
$\lambda\le0.$ In the same fashion one obtains that for $\lambda\le0,$
$\eta_m(\lambda)=0$ and $\zeta_m(\lambda)=0.$ Set $\chi_m(\lambda)=0$
for $\lambda\le0,$ and \eqref{6.10} will hold.

Let $\lambda^m_1 < \lambda^m_2 < \dots < \lambda^m_{P_m}$ be the
points in the interval $(0, d^2]$ where, according to lemma
\ref{lm6.2}, $M_0^m(\lambda)-M_{\infty}^m(\lambda)=0.$ We can use
lemma \ref{lm6.2} and \eqref{6.11} to deduce that $\xi_m,$ $\eta_m$ and
$\zeta_m$ are constant on each of the intervals making up $(0,
d^2]\backslash\{\lambda^m_1, \lambda^m_2, \dots, \lambda^m_{P_m}\}.$
Set $\chi_m$ to be constant on each of these intervals.  At each
of the points $\lambda^m_1, \lambda^m_2, \dots, \lambda^m_{P_m}$ the
functions $\xi_m,$ $\eta_m$ and $\zeta_m$ could have a jump. To show
\eqref{6.10} on $(0, d^2]$ we need to find a relationship between the
jumps of these functions. We will return to this shortly.

The remaining case, $\lambda > d^2,$ is treated similarly. Equation
\eqref{6.9} in the proof of lemma \ref{lm6.2} gives an expression for
$M_0^m(\lambda)-M_{\infty}^m(\lambda)$ on this interval (with the notation from
\eqref{4.11}).
By using the fact that the quantities $a_m,\, b_m$ and $j_m$ together with
their derivatives are real, we can calculate
\begin{equation*}
\begin{split}
  \mathfrak{Im}\, M(\lambda)
   & =  \mathfrak{Im}\Bigg\{\frac{j_m(a_m+ib_m)}{(j_m'a_m-a_m'j_m)+i(j_m'b_m-b_m'j_m)}\Bigg\}\\
   & =  \frac{j_m^2(b_m'a_m-a_m'b_m)}{(j_m'a_m-a_m'j_m)^2+(j_m'b_m-b_m'j_m)^2}.
\end{split}
\end{equation*}
Use \eqref{4.12} to simplify the numerator of this fraction.
Apply \eqref{6.11}. We get that for any $d^2 < \lambda_0 < \lambda_1,$
$$
  \xi_m(\lambda_1)-\xi_m(\lambda_0)=\frac{\pi}{2}\int\limits_{\lambda_0}^{\lambda_1}\!
  \frac{j_m(R,\lambda)^2}{c_m(\lambda)^2+d_m(\lambda)^2}\,d\lambda,
$$
where $c_m(\lambda)$ and $d_m(\lambda)$ are defined by \eqref{4.13a} and \eqref{4.13b}.
So we have
$$
  d\xi_m(\lambda)=\frac{\pi}{2}\frac{j_m(R,\lambda)^2\,d\lambda}{c_m(\lambda)^2+d_m(\lambda)^2}.
$$

We want \eqref{6.10} to hold. Define $\chi_m(\lambda)$ for $\lambda>d^2$ such that
$$
  d\chi_m(\lambda)=\frac{\pi}{2}\frac{d\lambda}{c_m(\lambda)^2+d_m(\lambda)^2},
$$
then the first of the three identities \eqref{6.10} is valid. It is easy
to repeat the same calculation for $\eta_m$ and $\zeta_m$ and show that
for $\lambda>d^2$ the other two identities in \eqref{6.10} hold.

Now we will return to what is the longest part of the proof, the study
of what happens at the points $\lambda\in (0, d^2]$ where
$M_0^m(\lambda)-M_{\infty}^m(\lambda)=0.$ Let $\lambda_0$ be such a point.  An immediate
observation is that $M_0^m(\lambda_0)$ is finite and $j_m(R,\lambda_0)\ne
0.$ Indeed, it was proved in lemma \ref{lm6.2} that $M_{\infty}^m(\lambda)$ is
finite for $\lambda$ real.  So, $M_0^m(\lambda_0)$ which equals $M_{\infty}^m(\lambda_0)$ is
finite.  Then, since $M_0^m(\lambda_0)=-j_m'(R,\lambda_0)/j_m(R,\lambda_0)$ and because
$j_m'(R,\lambda_0)$ cannot become zero simultaneously with $j_m(R,\lambda_0)$
($j_m(r,l)$ is a solution of \eqref{2.5}), we deduce
$j_m(R,\lambda_0)\ne0.$

With this observation in hand, in order to prove that \eqref{6.10}
holds at $\lambda_0$ one needs to show that
\begin{equation}\label{6.12}
  d\eta_m(\lambda_0)=M_0^m(\lambda_0)\,d\xi_m(\lambda_0), \ \ d\zeta_m(\lambda_0)=M_0^m(\lambda_0)^2\,d\xi_m(\lambda_0),
\end{equation}
and then define $d\chi_m(\lambda_0)=d\eta_m(\lambda_0)/j_m(R,\lambda_0)^2.$

Let $r_0$ be the jump of $\xi_m$ at $\lambda_0.$ Recall, $\xi_m$ was
defined by \eqref{5.9a}, so,
$$
  r_0=\lim\limits_{\varepsilon\to0^+}\lim\limits_{\delta\to0^+}
  \int\limits_{\lambda_0-\varepsilon}^{\lambda_0+\varepsilon}\!
  -\mathfrak{Im}\frac{1}{M_0^m(s+i\delta)-M_{\infty}^m(s+i\delta)}\,d s.
$$
By using \eqref{5.9b}, the first equality in \eqref{6.12} can be written
as
\begin{equation*}
  \lim\limits_{\varepsilon\to0^+}\lim\limits_{\delta\to0^+}
  \int\limits_{\lambda_0-\varepsilon}^{\lambda_0+\varepsilon}\!
  -\mathfrak{Im}\frac{M_0^m(s+i\delta)}{M_0^m(s+i\delta)-M_{\infty}^m(s+i\delta)}\,d s
  =M_0^m(\lambda_0)\,r_0.
\end{equation*}
$M_0^m(l)$ is analytic around $\lambda_0.$ Therefore, in a neighborhood of
$\lambda_0$
$$
  M_0^m(s+i\delta)=M_0^m(\lambda_0)+\{M_0^m(s)-M_0^m(\lambda_0)\}+i\delta H(s+i\delta),
$$
with $H$ a continuous function around $\lambda_0.$ Substitute this above.
By using the fact that $M_0^m(\lambda_0)$ is real and the definition of $r_0,$
we get the equivalent equality
\begin{multline}\label{6.13}
  M_0^m(\lambda_0)\,r_0\\
 + \lim\limits_{\varepsilon\to0^+}\lim\limits_{\delta\to0^+}
  \int\limits_{\lambda_0-\varepsilon}^{\lambda_0+\varepsilon}\!
  -\mathfrak{Im}\frac{\{M_0^m(s)-M_0^m(\lambda_0)\}+i\delta H(s+i\delta)}{M_0^m(s+i\delta)-M_{\infty}^m(s+i\delta)}\,d s\\
  =M_0^m(\lambda_0)\,r_0,
\end{multline}
so we have to prove that the limit of the integral on the left-hand
side of \eqref{6.13} is zero. To do this, we will need additional
information about $M_0^m-M_{\infty}^m.$

Let $l=s+i\delta$ $(\delta>0).$ Recall that formulas \eqref{5.6}
and \eqref{5.8} hold, where $\psi^m_0$ and $\psi^m_\infty$ are defined in
\eqref{5.5} and \eqref{5.7} and an expression for them is given by
\eqref{6.2}. Add equalities \eqref{5.6} and \eqref{5.8}. Obtain
$$
  \frac{1}{\delta}\mathfrak{Im}\{M_0^m(l)-M_{\infty}^m(l)\}=
  \int\limits_{0}^{R}\! \bigg|\frac{j_m(r, l)}{j_m(R, l)}\bigg|^2\,dr +
  \int\limits_{R}^{\infty}\! \bigg|\frac{w_m(r, l)}{w_m(R,l)}\bigg|^2\,dr.
$$
This equality gives us two things. First, that
$\mathfrak{Im}\{M_0^m(l)-M_{\infty}^m(l)\}>0,$ therefore,
\begin{equation}\label{6.14}
-\mathfrak{Im}\frac{1}{M_0^m(l)-M_{\infty}^m(l)}>0.
\end{equation}
Second, for $l$ close to $\lambda_0,$ the quantity
$\delta^{-1}\mathfrak{Im}\{M_0^m(l)-M_{\infty}^m(l)\}$ is bounded from below by
a strictly positive number, say $\omega_\infty^{-1}.$ Then
$\delta^{-1}|M_0^m(l)-M_{\infty}^m(l)|$ is bounded below by the same
number and therefore,
\begin{equation}\label{6.15}
\frac{\delta}{|M_0^m(l)-M_{\infty}^m(l)|}\le \omega_\infty.
\end{equation}

The integral in \eqref{6.13} can be written as
\begin{multline}\label{6.16}
  \int\limits_{\lambda_0-\varepsilon}^{\lambda_0+\varepsilon}\!
  -\mathfrak{Im}\frac{M_0^m(s)-M_0^m(\lambda_0)}{M_0^m(s+i\delta)-M_{\infty}^m(s+i\delta)}\,d s\\
+  \int\limits_{\lambda_0-\varepsilon}^{\lambda_0+\varepsilon}\!
  -\mathfrak{Im}\frac{i\delta H(s+i\delta)}{M_0^m(s+i\delta)-M_{\infty}^m(s+i\delta)}\,d s.
\end{multline}
Denote $c_\varepsilon=\sup_{|s-\lambda_0|\le\varepsilon}|M_0^m(s)-M_0^m(\lambda_0)|$. Since $M_0^m$ is continuous
in a neighborhood of $\lambda_0,$ we will have $c_\varepsilon\to 0$ as $\varepsilon\to 0.$ Let us estimate the
first integral from \eqref{6.16}.
\begin{equation*}
\begin{split}
  \bigg|\int\limits_{\lambda_0-\varepsilon}^{\lambda_0+\varepsilon}\!
  -\mathfrak{Im} \frac{M_0^m(s)-M_0^m(\lambda_0)}{M_0^m(s+i\delta)-M_{\infty}^m(s+i\delta)}\bigg|\,d s
   \le  \int\limits_{\lambda_0-\varepsilon}^{\lambda_0+\varepsilon}\!\bigg|
  \mathfrak{Im}\frac{M_0^m(s)-M_0^m(\lambda_0)}{M_0^m(s+i\delta)-M_{\infty}^m(s+i\delta)}\bigg|\,d s \\
  = \int\limits_{\lambda_0-\varepsilon}^{\lambda_0+\varepsilon}\!|M_0^m(s)-M_0^m(\lambda_0)|\,\bigg|
  \mathfrak{Im}\frac{1}{M_0^m(s+i\delta)-M_{\infty}^m(s+i\delta)}\bigg|\,d s\\
  \le c_\varepsilon  \int\limits_{\lambda_0-\varepsilon}^{\lambda_0+\varepsilon}\!\bigg|
  \mathfrak{Im}\frac{1}{M_0^m(s+i\delta)-M_{\infty}^m(s+i\delta)}\bigg|\,d s\\
  = c_\varepsilon \int\limits_{\lambda_0-\varepsilon}^{\lambda_0+\varepsilon}\!
  -\mathfrak{Im}\frac{1}{M_0^m(s+i\delta)-M_{\infty}^m(s+i\delta)}\,d s\\
  = c_\varepsilon [\xi_m(\lambda_0+\varepsilon)-\xi_m(\lambda_0-\varepsilon)].
 \end{split}
 \end{equation*}
In deriving this we used the fact that $M_0^m(s)-M_0^m(\lambda_0)$ is
real, together with \eqref{6.14} and \eqref{5.9a}. Clearly as
$\varepsilon\to0,$ the integral goes to zero.

Now estimate the second integral in the sum \eqref{6.16}. Let $H_\infty$
be an upper bound of $|H(s+i\delta)|$ for $l=s+i\delta$ in a
neighborhood of $\lambda_0.$ Inequality \eqref{6.15} implies that as  $\varepsilon\to0,$
$$
  \bigg|\int\limits_{\lambda_0-\varepsilon}^{\lambda_0+\varepsilon}\!
  -\mathfrak{Im}\frac{i\delta H(s+i\delta)}{M_0^m(s+i\delta)-M_{\infty}^m(s+i\delta)}\,d s\bigg|\le
  \int\limits_{\lambda_0-\varepsilon}^{\lambda_0+\varepsilon}\!H_\infty \omega_\infty\,d s
  =2\varepsilon H_\infty \omega_\infty\to 0
$$
This proves the first equality in \eqref{6.12}.

Let us prove the second equality in \eqref{6.12}. Recall that $\zeta_m$ is
given by \eqref{5.9c}. The above approach does not apply immediately,
since unlike $M_0^m(l)$ (the numerator in \eqref{5.9b}), the
function $M_0^m(l) M_{\infty}^m(l)$ (the numerator in \eqref{5.9c})
will not be analytic around $\lambda_0$ if $\lambda_0=d^2$ (since as seen from lemma \ref{lm3.4},
$w_m(r, l)$ is not defined in a neighborhood of  $l=d^2$).
The idea is then to use the equality
$$
  \frac{x y}{x-y}=\frac{x^2}{x-y}-x,
$$
to write the jump of $\zeta_m$ at $\lambda_0$ as
$$
  \lim\limits_{\varepsilon\to0^+}\lim\limits_{\delta\to0^+}\!\!
  \int\limits_{\lambda_0-\varepsilon}^{\lambda_0+\varepsilon}\!\!
  -\mathfrak{Im}\frac{M_0^m(s+i\delta)^2}{M_0^m(s+i\delta)-M_{\infty}^m(s+i\delta)}\,d s
  +\lim\limits_{\varepsilon\to0^+}\lim\limits_{\delta\to0^+}\!\!
  \int\limits_{\lambda_0-\varepsilon}^{\lambda_0+\varepsilon}\!\!
  \mathfrak{Im}\, M_0^m(s+i\delta)\,d s.
$$
The limit of the second integral in the sum will be zero, since as
argued above, $M_0^m$ will be finite at $\lambda_0$ (and therefore,
around $\lambda_0$) and thus the quantity inside the integral is
bounded.  For the first integral in the sum we can proceed in the same
way we calculated the jump of $\eta_m.$ This finishes the proof of
\eqref{6.10}.

Finally, we need to justify the claim that $\chi_m$ is a
non-decreasing function. From \eqref{6.10} we have
$$
  d\xi_m(\lambda)+d\zeta_m(\lambda)=\{j_m(R,\lambda)^2+j_m'(R,\lambda)^2\}d\chi_m(\lambda).
$$
The left-hand side of this is non-negative measure, since by theory
$\xi_m$ and $\zeta_m$ are non-decreasing. The number
$j_m(R,\lambda)^2+j_m'(R,\lambda)^2$ is strictly positive, as
$j_m(R,\lambda)$ and $j_m'(R,\lambda)$ cannot be both zero,
$j_m(r,\lambda)$ being a non-zero solution to the second order
differential equation \eqref{2.5}. Then we get that $d\chi_m(\lambda)$
is a non-negative measure, which shows that $\chi_m$ is
non-decreasing.
\end{proof}

\begin{corollary}\label{cor6.4}
Let $\lambda\in(0, d^2]$ be a discontinuity point for $\chi_m.$ Then, if $\lambda<d^2,$
the following hold:
\begin{subequations}
\begin{equation}\label{6.17a}
\displaystyle  \frac{j_m'(R, \lambda)}{j_m(R, \lambda)}=\frac{k_m'(R, \lambda)}{k_m(R, \lambda)},
\end{equation}
and
\begin{equation}\label{6.17b}
  j_m(r,\lambda)=\displaystyle \frac{j_m(R,\lambda)}{k_m(R,\lambda)}k_m(r,\lambda),\ \ r\ge R.
\end{equation}
\end{subequations}
While for $\lambda=d^2,$
\begin{subequations}
\begin{equation}\label{6.18a}
\displaystyle\frac{j_m'(R, \lambda)}{j_m(R, \lambda)}=-\displaystyle\frac{|m|-1/2}{R},
\end{equation}
and
\begin{equation}\label{6.18b}
  j_m(r,\lambda)=\displaystyle\frac{j_m(R,\lambda)}{R^{1/2-|m|}}r^{1/2-|m|},\ \ r\ge R.
\end{equation}
\end{subequations}
In particular, for $\lambda\in(0, d^2]$ a discontinuity point of $\chi_m,$
$j_m(r,\lambda)$ decays exponentially as $r\to\infty$ if $\lambda<d^2,$ and $j_m(r,\lambda)\sim r^{1/2-|m|}$
as $r\to\infty$ if $\lambda=d^2.$
\end{corollary}
\begin{proof}
As it follows from the proof of Theorem \ref{th6.3}, for such a $\lambda$ we will
have $M_0^m(\lambda_0)-M_\infty^m(\lambda_0)=0.$ With the help \eqref{6.1}, \eqref{6.7a} and \eqref{6.8},
we deduce  \eqref{6.17a} and \eqref{6.18a}.  Notice that these equalities
are exactly the conditions \eqref{4.6} and \eqref{4.9}. Then  \eqref{6.17b} and \eqref{6.18b}
follow from \eqref{4.7} and \eqref{4.10}.
\end{proof}

\section{Computing the Transform}\label{sec7}
Denote by $L^2(\chi_m)$ the space of all functions $G:\mathbb
R\to\mathbb C$ such that
$$
\int_{-\infty}^{\infty}\!|G(\lambda)|^2\,d\chi_m(\lambda)<\infty,
$$
where $\chi_m$ is the non-decreasing function
defined in Theorem \ref{th6.3}.
\begin{theorem}\label{th7.1}
Let $g\in L^2(0,\infty).$  The integral
\begin{equation}\label{7.1}
  G_m(\lambda)=\int\limits_{0}^{\infty}\!j_m(r,\lambda)g(r) \,dr
\end{equation}
is convergent in $L^2(\chi_m)$, in the sense that there exists $G_m\in
L^2(\chi_m)$ such that $G_m^{cd}\to G_m$ in $L^2(\chi_m)$ as $c\to 0$ and $d\to\infty,$ where
\begin{equation}\label{7.2}
  G_m^{cd}(\lambda)=\int\limits_{c}^{d}\!j_m(r,\lambda)g(r) \,dr, \  \ 0 < c < d <\infty.
\end{equation}
The equality
\begin{equation}\label{7.3}
  g(r)= \frac{1}{\pi} \int\limits_{-\infty}^{\infty}\! j_m(r,\lambda)G_m(\lambda)\, d\chi_m(\lambda)
\end{equation}
holds, in the sense that $g^{\sigma\tau}\to g$ in $L^2(0,\infty)$ as
$\tau\to-\infty,\, \sigma\to\infty,$ where
\begin{equation}\label{7.4}
g^{\tau\sigma}(r)=\frac{1}{\pi}\int\limits_{\tau}^{\sigma}\!j_m(r,\lambda)G_m(\lambda)\,d\chi_m(\lambda),\ \
                    -\infty<\sigma<\tau<\infty.
\end{equation}
We have the Parseval identity
\begin{equation}\label{7.5}
  \int\limits_{0}^{\infty}\!|g(r)|^2\,dr=
  \frac{1}{\pi}\int\limits_{-\infty}^{\infty}\!|G_m(\lambda)|^2 \,d\chi_m(\lambda).
\end{equation}
\end{theorem}

\begin{proof} We will apply Theorem \ref{th5.1}. First note that if
$\Gamma=(\Gamma_1, \Gamma_2)$ with $\Gamma_1,\,\Gamma_2\,:\mathbb R\to \mathbb C,$ then
because of \eqref{6.10}, the norm $||\cdot||$ as defined in
\eqref{5.11} can be written in the form
\begin{equation}\label{7.6}
  ||\Gamma||^2=\int\limits_{-\infty}^{\infty}\!
 |j_m(R,\lambda)\Gamma_1(\lambda)-j_m'(R,\lambda)\Gamma_2(\lambda)|^2\,d\chi_m(\lambda).
\end{equation}

Recall the identities \eqref{6.3} to \eqref{6.5} for expressing
$\theta(r,\lambda)$ and $\varphi(r,\lambda)$ in terms of
$j_m(r,\lambda)$ and $y_m(r,\lambda)$.  Note that from \eqref{6.4} and
\eqref{6.5} we get the following equalities involving the coefficients
$\alpha, \beta, \gamma,$ and $\delta$
\begin{equation}\label{7.7}
  \alpha(\lambda)j_m(R,\lambda)-\gamma(\lambda)j_m'(R,\lambda)=1,
\end{equation}
and
\begin{equation}\label{7.8}
  \beta(\lambda)j_m(R,\lambda)-\delta(\lambda)j_m'(R,\lambda)=0.
\end{equation}
Then, for $0 < c< d < \infty,$ the functions $\Gamma^{cd}_1$ and
$\Gamma^{cd}_2$ defined by \eqref{5.12} become
$$
  \Gamma^{cd}_1(\lambda)=\int\limits_{c}^{d}\!\big\{\alpha(\lambda)j_m(r,\lambda)+\beta(\lambda)y_m(r, \lambda)\big\} g(r)\,dr,
$$
$$
  \Gamma^{cd}_2(\lambda)=\int\limits_{c}^{d}\!\big\{\gamma(\lambda)j_m(r,\lambda)+\delta(\lambda)y_m(r, \lambda)\big\}g(r)\,dr.
$$
We denoted $\Gamma^{cd}=(\Gamma^{cd}_1,\Gamma^{cd}_2).$ Let
\begin{equation}\label{7.9}
  \tilde \Gamma^{cd}_1(\lambda)=\int\limits_{c}^{d}\!\alpha(\lambda)j_m(r, \lambda) g(r)\,dr, \ \
  \tilde \Gamma^{cd}_2(\lambda)=\int\limits_{c}^{d}\!\gamma(\lambda)j_m(r, \lambda) g(r)\,dr.
\end{equation}
and set $\tilde \Gamma^{cd}=(\tilde \Gamma^{cd}_1,\tilde \Gamma^{cd}_2).$ We
have
\begin{equation}\label{7.10}
  ||\Gamma^{cd}-\tilde \Gamma^{cd}||=0.
\end{equation}
That follows from \eqref{7.6} and \eqref{7.8}.  Theorem \ref{th5.1}
guarantees the existence of $\Gamma=(\Gamma_1,\Gamma_2)\in L^2(\rho)$ such that
$$
  ||\Gamma-\Gamma^{cd}||\to 0 \mbox{ as } c\to 0, \, d\to \infty.
$$
Then, equality \eqref{7.10} says that we have
$||\Gamma-\tilde \Gamma^{cd}||\to 0 \mbox{ as } c\to 0, \, d\to \infty.$
But, according to \eqref{7.6},
\begin{multline}\label{7.11}
  ||\Gamma-\tilde \Gamma^{cd}||^2
  =\int\limits_{-\infty}^{\infty}\!
  |\{j_m(R,\lambda)\Gamma_1(\lambda)-j_m'(R,\lambda)\Gamma_2(\lambda)\}-\\
  \{j_m(R,\lambda)\tilde \Gamma^{cd}_1(\lambda)-j_m'(R,\lambda)\tilde \Gamma^{cd}_2(\lambda)\}|^2\,d\chi_m.
\end{multline}
Note that
$$
  j_m(R,\lambda)\tilde \Gamma^{cd}_1(\lambda)-j_m'(R,\lambda)\tilde \Gamma^{cd}_2(\lambda)
  =\int\limits_{c}^{d}\! j_m(r, \lambda) g(r)\,dr.
$$
That follows from \eqref{7.7} and \eqref{7.9}. Thus, if we denote
\begin{equation}\label{7.12}
  G_m(\lambda)=j_m(R,\lambda)\Gamma_1(\lambda)-j_m'(R,\lambda)\Gamma_2(\lambda),\, \lambda\in\mathbb R,
\end{equation}
we obtain from \eqref{7.11} that
$$
  \int\limits_{-\infty}^{\infty}\!|G_m(\lambda)-G_m^{cd}(\lambda)|^2  \,d\chi_m\to 0
  \mbox{ as } c\to0,\, d\to\infty
$$
($G_m^{cd}$ was defined by \eqref{7.2}). This shows the first part of
Theorem \ref{th7.1}.

Next we need to show that  representation \eqref{7.3} holds.  It
suffices to prove that $g^{\tau\sigma}$ as defined by \eqref{5.13} in
Theorem \ref{th5.1} is the same as $g^{\tau\sigma}$ defined in
\eqref{7.4}. And they are. To check this one needs to start with
$g^{\tau\sigma}$ as given in \eqref{5.13}, substitute $\theta(r,\lambda)$
and $\varphi(r,\lambda)$ from \eqref{6.3}, use \eqref{6.10} to
express $\xi_m,\,\eta_m$ and $\zeta_m$ in terms of $\chi_m,$ use the
equalities \eqref{7.7} and \eqref{7.8}, and finally use definition
\eqref{7.12} for $G_m(\lambda).$

Lastly, the Parseval identity \eqref{7.5} follows from \eqref{5.14},
\eqref{7.6} and \eqref{7.12}. The theorem is proved.
\end{proof}

\begin{theorem}\label{th7.2}
Let $g\in L^2(0,\infty).$ Let $\chi_m$ be the non-decreasing function defined in theorem
{\rm \ref{th6.3}}.  Let $0< \lambda^m_1< \dots < \lambda^m_{P_m}\le d^2$
$(P_m\ge 0)$ be the points where $\chi_m$ is discontinuous. Let
$r^m_1,\dots,r^m_{P_m}$ be the corresponding jumps. Let
$a_m(r,\lambda)$ and $b_m(r,\lambda)$ be the functions defined by
{\rm \eqref{4.11}}, and $c_m(\lambda)$ and $d_m(\lambda)$ be defined by
{\rm \eqref{4.13a}} and {\rm \eqref{4.13b}}.  Then,
\begin{equation}\label{7.13}
  r^m_k=\pi \bigg\{\int\limits_{0}^{\infty}\!j_m(r,\lambda^m_k)^2\,dr\bigg\}^{-1},\, k=1,\dots, P_m,
\end{equation}
\begin{equation}\label{7.14}
  d\chi_m(\lambda)=\frac{\pi}{2}\frac{d\lambda}{c_m(\lambda)^2+d_m(\lambda)^2}, \ \ \lambda\in (d^2, \infty).
\end{equation}
We have the representation
\begin{equation}\label{7.15}
  g(r)=\frac{1}{\pi}\sum\limits_{k=1}^{P_m} r^m_k j_m(r,\lambda^m_k)G_m(\lambda^m_k)+
  \frac{1}{2}\int\limits_{d^2}^{\infty}\!\frac{j_m(r,\lambda)G_m(\lambda)}
                                                        {c_m(\lambda)^2+d_m(\lambda)^2}\,d\lambda.
\end{equation}
\end{theorem}
\begin{proof}  In theorem \ref{th6.3} we proved the existence of the
function $\chi_m$ and along the way we found its continuous part, that is,
equality \eqref{7.14}.  We need to find its discrete part, that is,
the value of all the jumps of $\chi_m.$ Then \eqref{7.15} will
follow by applying \eqref{7.3}.

Let us notice the following observation.  If
$\lambda_1\ne\lambda_2,$ and $w_1(r)$ and $w_2(r)$ satisfy \eqref{2.5}
with $\lambda=\lambda_1,$ and $\lambda=\lambda_2$ respectively, then
for any $0< c < d <\infty$
\begin{equation}\label{7.16}
  \int\limits_{c}^{d}\!w_{1}(r)w_2(r)\,dr
  =-(\lambda_1-\lambda_2)^{-1}[w'_1(r)w_2(r)-w_1(r)w'_2(r)]\bigg|^c_d.
\end{equation}
Indeed, we can write
$$
  \int\limits_{c}^{d}\!w''_1(r)w_2(r)\,dr=-\int\limits_{c}^{d}\!
  \bigg\{\lambda_1 - q(r) - \frac{m^2-1/4}{r^2}\bigg\}w_1(r)w_2(r)\,dr,
$$
and
$$
  \int\limits_{c}^{d}\!w_1(r)w''_2(r)\,dr=-\int\limits_{c}^{d}\!
  w_1(r)\bigg\{\lambda_2 - q(r) - \frac{m^2-1/4}{r^2}\bigg\}w_2(r)\,dr.
$$
If we integrate by parts the left-hand sides of these two equalities,
and then subtract from first the second, we get exactly \eqref{7.16}.

Before we prove \eqref{7.13}, recall the behavior of $j_m(r,\lambda)$ as
$r\to\infty.$ For $\lambda\le d^2$ a discontinuity point of $\chi_m$
it is described in corollary \ref{cor6.4}, while for $\lambda>d^2$
see the discussion in section 4.

Let $\lambda_0\le d^2$ be one of the discontinuity points of $\chi_m.$
Let $r_0$ be the corresponding jump.  We will consider two cases: when
$j_m(r,\lambda_0)$ is square integrable, and when it is not. The first
case splits into two subcases: we can have either $\lambda_0< d^2,$ or
$\lambda_0=d^2$ with $|m|\ge 2.$ The second case happens for
$\lambda_0=d^2$ and $|m|\in\{0,1\}.$ Let us start with the first case.

Denote $g(r)=j_m(r,\lambda_0).$ Since $g(r)$ is square integrable, we can
apply theorem \ref{th7.1} for this particular function.  We will show
that the corresponding $G_m(\lambda)$ as defined by \eqref{7.1} is such that
\begin{equation}\label{7.17}
G_m(\lambda)=\left\{
\begin{array}{cl}
\int_{0}^{\infty}\!j_m(r,\lambda_0)^2\,dr, & \mbox{ if } \lambda=\lambda_0,\\
0,                                   & \mbox{ if } \lambda\ne\lambda_0.
\end{array}
\right.
\end{equation}
for all $\lambda$ such that $d\chi_m(\lambda)\ne0.$ Then \eqref{7.13} will follow
promptly; one needs to apply the Parseval identity \eqref{7.5} and notice
that since $r_0$ is the jump of $\chi_m$ at $\lambda_0,$ then
$d\chi_m(\lambda_0)=r_0\delta(\lambda-\lambda_0)$ with $\delta$ being Dirac's
function.

Consider first the subcase $\lambda_0< d^2.$ That \eqref{7.17} is true
for $\lambda=\lambda_0$ follows immediately from \eqref{7.1}. Assume now
$\lambda\ne\lambda_0.$ Let us compute $G_m^{cd}(\lambda)$ for $0< c < d
< \infty$ as defined in \eqref{7.2}.  Use \eqref{7.16} with
$w_1(r)=j_m(r,\lambda_0),$ $w_2(r)=j_m(r,\lambda).$ Put $c\to0$ and
$d\to\infty.$ From \eqref{3.1} we deduce that
$$
  j_m'(c,\lambda_0)j_m(c,\lambda)-j_m(c,\lambda_0)j_m'(c,\lambda) \to 0 \mbox{ as } c\to0.
$$
Also, both $j_m'(d,\lambda_0)j_m(d,\lambda)$ and
$j_m(d,\lambda_0)j_m'(d,\lambda)$ go to zero as $d\to\infty$ if
$d\chi_m(\lambda)\ne 0,$ since on one hand, $j_m(r,\lambda_0)$ and its
derivative decrease exponentially, and on the other hand,
$j_m'(r,\lambda)$ and its derivative either decrease exponentially for
$\lambda<d^2$, or behave like a power of $r$ for $\lambda=d^2,$ or are
bounded for $\lambda > d^2.$ In any case we get that
$G_m^{cd}(\lambda)\to 0$ as $c\to0$ and $d\to\infty,$ so
$G_m(\lambda)=0.$

The subcase $\lambda_0=d^2$ with $|m|\ge 2$ follows in the same way.
The statement that for $\lambda\ne\lambda_0$ and $d\chi_m(\lambda)\ne
0$ both $j_m'(d,\lambda_0)j_m(d,\lambda)$ and
$j_m(d,\lambda_0)j_m'(d,\lambda)$ go to zero as $d\to\infty$ is argued
in a little different way. We have that $j_m(r,\lambda_0)$ and its
derivative decay as a negative power of $r$ for $r\to\infty,$ while
$j_m(r,\lambda)$ and its derivative either decay exponentially for
$\lambda< \lambda_0,$ or stay bounded for $\lambda> \lambda_0.$ But
the conclusion is the same, $G_m^{cd}(\lambda)\to 0$ as $c\to0$ and
$d\to\infty,$ and thus \eqref{7.17} holds in this case too.

Now consider the second case, $\lambda_0=d^2$ and $|m|\in\{0,1\},$ when,
as we remarked above, $j_m(r, \lambda_0)$ is not square
integrable. For $0 < c < d < \infty$ define
$$
  g(r)=\left\{
  \begin{array}{cl}
  j_m(r, \lambda_0), & \mbox{ if } c < r <  d,\\
  0,                                   & \mbox{ otherwise.}
  \end{array}
  \right.
$$
This function \emph{will be} square integrable. Let $G_m(\lambda)$ be
the corresponding transform of $g(r)$ as defined by \eqref{7.1}. Apply
the Parseval identity \eqref{7.5}. Get
\begin{multline*}
  \int\limits_{c}^{d}\!|j_m(r,\lambda_0)|^2\,dr
  =\frac{1}{\pi}\int\limits_{-\infty}^{\infty}\!|G_m(\lambda)|^2 \,d\chi_m(\lambda)\\
  \ge \frac{1}{\pi}\int\limits_{\{\lambda_0\}}|G_m(\lambda)|^2 \,d\chi_m(\lambda)
 =\frac{1}{\pi}G_m(\lambda_0)^2 r_0.
\end{multline*}
From \eqref{7.1} we obtain that
$$
  G_m(\lambda_0)=\int_{c}^{d}\!|j_m(r,\lambda_0)|^2\,dr.
$$
Then we can write
$$
  \pi \bigg\{ \int\limits_{c}^{d}\!|j_m(r,\lambda_0)|^2\,dr\bigg\}^{-1}\ge r_0.
$$
By putting $c\to0,$ $d\to\infty$ and noticing that $r_0\ge0,$ being a
jump of the non-decreasing function $\chi_m,$ we deduce $r_0=0.$
\end{proof}

With this theorem, corollary \ref{cor6.4} and the classification
obtained in section 4, we can characterize completely the points of
discontinuity of $\chi_m.$
\begin{corollary}\label{cor7.3} Let $\lambda\in \mathbb R.$
Then $\lambda$ is a discontinuity point of  $\chi_m$ if and
only if $j_m(\cdot, \lambda)\in L^2(0,\infty),$ and if and only if
$0< \lambda<d^2$ and {\rm \eqref{6.17a}} holds, or $\lambda=d^2$ and {\rm \eqref{6.18a}} holds.
\end{corollary}

\section{Finding Green's Function}\label{sec8}
In this section we will show that under certain conditions, the
solution of the Helmholtz equation \eqref{1.1}, which in the
cylindrical coordinate system $(r,\vartheta, z)$ is written as
\eqref{2.1}, is unique.  We will find a representation for it in terms
of the source $f(r,\vartheta, z),$ the eigenfunctions $j_m(r,\lambda)$ of
equation \eqref{2.5} satisfying lemma \ref{lm3.1} and the measure
$d\chi_m(\lambda)$ defined in theorem \ref{th6.3}.  Before proving
this result we will need one lemma.

\begin{lemma}\label{lm8.1} Let $u\in C^1(\mathbb R^3).$ Then $u$ can be written as
\begin{equation}\label{8.1}
   u(r,\vartheta, z)=\sum\limits_{m\in \mathbb Z} e^{im\vartheta}u_{m}(r,z).
\end{equation}
For each $z\in \mathbb R,$ the function $r\to u_m(r,z)$ is in $C^1
[0, \infty)$ and
\begin{equation}\label{8.2}
  \lim\limits_{r\to0}  \bigg[j_m(r,\lambda)\frac{\partial \{\sqrt{r}u_m(r,z)\}}{\partial r}-
  \frac{\partial j_m(r,\lambda)}{\partial r}\{\sqrt{r}u_m(r,z)\}\bigg]=0.
\end{equation}
\end{lemma}
\begin{proof}
Equality \eqref{8.1} is nothing but the Fourier series of the function $\vartheta\to u(r,\vartheta, z).$
The smoothness  of the obtained $u_m$ follows from the formula for the Fourier coefficients,
$$
  u_{m}(r,z)=\frac{1}{2\pi}\int\limits_{0}^{2\pi}\!u(r,t, z)e^{-imt}\,dt.
$$

Let us prove \eqref{8.2}. Write it as
$$
  \lim\limits_{r\to0}  \bigg[
  \sqrt{r}j_m(r,\lambda)\, \frac{\partial u_m(r,z)}{\partial r} + u_m(r,z)
  \bigg(\frac{j_m(r,\lambda)}{2\sqrt{r}}-\sqrt{r} \frac{\partial j_m(r,\lambda)}{\partial r}\bigg)\bigg]=0.
$$
The first term in the sum clearly goes to zero as $r\to0.$ For the second term,
by applying lemma \ref{lm3.1} to the expression in parentheses we get
\begin{multline*}
  \lim\limits_{r\to0}
  \bigg(\frac{j_m(r,\lambda)}{2\sqrt{r}}-\sqrt{r} \frac{\partial j_m(r,\lambda)}{\partial r}\bigg)\\
  =\lim\limits_{r\to0}
  \bigg( \frac{1}{2}r^{|m|} - \bigg\{ |m|+\frac{1}{2} \bigg\} r^{|m|}\bigg)
  = -\lim\limits_{r\to0} |m| r^{|m|} =0,
\end{multline*}
for all $m\in \mathbb Z.$
\end{proof}

These will be the conditions on the solution $u$ for \eqref{2.1} which
will guarantee its uniqueness. First, we will assume that the source
$f$ is continuous and with compact support. Second, we will impose the
condition that $u \in C^1 (\mathbb R^3).$ Third,
suppose that for all $m\in\mathbb Z,$ $z\in\mathbb R$ the following
equality holds
\begin{equation}\label{8.3}
  \lim\limits_{r\to\infty}  \bigg[j_m(r,\lambda)\frac{\partial \{\sqrt{r}u_m(r,z)\}}{\partial r}-
  \frac{\partial j_m(r,\lambda)}{\partial r}\{\sqrt{r}u_m(r,z)\}\bigg]=0,
\end{equation}
with the functions $u_m(r, z)$ defined by \eqref{8.1}.

Denote by $U_m(\lambda, z)$ the transform of the function $r\to
\sqrt{r}u_m(r,z)$
given by \eqref{7.1},
\begin{equation}\label{8.4}
  U_m(\lambda, z)=\int\limits_{0}^{\infty}\!j_m(\rho,\lambda) \sqrt{r}u_m(r,z)\,d\rho,
\end{equation}
The fourth requirement is the \emph{radiation condition}
\begin{equation}\label{8.5}
\begin{cases} \displaystyle \lim_{|z| \to \infty} \left[ \frac{\partial U_m (\lambda,z)}{\partial |z|} - i \sqrt{k^2 n_0^2 - \lambda} \;
U_m (\lambda,z)
\right] = 0,  &\textmd{for } \lambda \le k^2 n_0^2 \mbox{ with } d\chi_m (\lambda)\!\ne\! 0, \\  & \\
\displaystyle \ \lim_{|z| \to \infty} U_m (\lambda,z) = 0, & \textmd{for } \lambda > k^2 n_0^2.
\end{cases}
\end{equation}

These conditions are physically motivated. First condition says that the
source must be finite. Equation \eqref{8.3} signifies a fast decay of
the electromagnetic field intensity as $r\to\infty$. And the radiation
condition \eqref{8.5} means that the energy going to $z=\infty$ can be
separated in two parts. First part is oscillatory, and it goes to
infinity, and \emph{cannot not} come from infinity, while the second
part is rapidly decaying.

\begin{theorem}\label{th8.2}
With the above assumptions, the solution of \eqref{2.1} can be
represented as
\begin{equation}\label{8.6}
  u(r,\vartheta, z)=\int\limits_{-\infty}^{\infty}\!\int\limits_{0}^{\infty}\!\int\limits_{0}^{2\pi}\!
  G(r, \rho; \vartheta, t; z, \zeta) f(\rho,t, \zeta)\,\rho \,dt\,d\rho\,d\zeta,
\end{equation}
where
\begin{multline}\label{8.7}
  G(r, \rho; \vartheta, t; z, \zeta)=\\
  \frac{1}{2\pi^2}\frac{1}{\sqrt{r\rho}}
  \sum\limits_{m\in\mathbb Z}\ \int\limits_{-\infty}^{+\infty}\!
  \frac{e^{i |z-\zeta|\sqrt{k^2n_0^2-\lambda}}}{2i\sqrt{k^2n_0^2-\lambda}}
  e^{i m (\vartheta-t)}\,  j_m(\rho,\lambda)\, j_{m}(r,\lambda)\,d\chi_m(\lambda),\\
0 < r, \rho;\, 0 \le \vartheta, t\le 2\pi;\, z, \zeta\in \mathbb R,
\end{multline}
and $\chi_m$ is the non-decreasing function defined in Theorem \ref{th6.3}.
\end{theorem}
\begin{proof}
The function $f(\rho,\vartheta, z)$ can be decomposed as
\begin{equation}\label{8.8}
  f(\rho,\vartheta, z)=\sum\limits_{m\in\mathbb Z} e^{im\vartheta}f_{m}(\rho,z),
\end{equation}
with
\begin{equation}\label{8.9}
  f_{m}(\rho,z)=\frac{1}{2\pi}\int\limits_{0}^{2\pi}\!f(\rho,t, z)e^{-imt}\,dt.
\end{equation}
Look for $u(\rho,\vartheta, z)$ in the form \eqref{8.1}.  By
plugging \eqref{8.1} and \eqref{8.8} into \eqref{2.1} we deduce that
for each $m,$ $u_m(\rho,z)$ needs to satisfy the equation
\begin{equation}\label{8.10}
\frac{\partial^2 u_m}{\partial z^2}+\frac{1}{\rho}\frac{\partial}{\partial \rho}
\bigg(\rho \frac{\partial u_m}{\partial \rho}\bigg)+\bigg(k^2 n(\rho)^2-\frac{m^2}{\rho^2}\bigg)u_m=f_m,
\end{equation}
for all $m\in\mathbb Z.$
Let $\lambda\in\mathbb R$ be such that $d\chi_m(\lambda)\ne 0.$ Let
$U_m(\lambda, z)$ be the transform of $u_m(\rho, z)$ given by
\eqref{8.4}, and let
\begin{equation}\label{8.11}
   F_m(\lambda, z)=\int\limits_{0}^{\infty}\!j_m(\rho,\lambda) \sqrt{\rho}f_m(\rho,z)\,d\rho
\end{equation}
be the transform of $f_m.$ Multiply \eqref{8.10} on both sides by
$\sqrt{\rho}j_m(\rho,\lambda)$ and integrate from $0$ to $\infty.$
Obtain
\begin{multline*}
  \frac{\partial^2 U_m}{\partial z^2} +
  \int\limits_{0}^{\infty}\!j_m(\rho,\lambda) \frac{1}{\sqrt{\rho}}\frac{\partial}{\partial \rho}
  \bigg(\rho \frac{\partial u_m(\rho,z)}{\partial \rho}\bigg)\,d\rho+\\
  \int\limits_{0}^{\infty}\!j_m(\rho,\lambda)
  \bigg(k^2 n(\rho)^2-\frac{m^2}{\rho^2}\bigg)\sqrt{\rho}u_m(\rho,z)\,d\rho
  =  F_m.
\end{multline*}
Use integration by parts twice for the first integral in the above equation.  By
applying lemma \ref{lm8.1}, and equality
\eqref{8.3}  we get
$$
  \frac{\partial^2 U_m}{\partial z^2} +
  \int\limits_{0}^{\infty}\!\bigg\{\frac{\partial^2 j_m(\rho,\lambda)}{\partial \rho}+\bigg(k^2 n(\rho)^2 -
  \frac{m^2-1/4}{\rho^2}\bigg)j_m(\rho,\lambda)\bigg\} \sqrt{\rho}u_m(\rho,\lambda) \,d\rho
  =  F_m.
$$
Recall that $j_m(\rho,\lambda)$ satisfies \eqref{2.5} with $q(r)$
given by \eqref{2.3}. Then,
\begin{equation}\label{8.12}
  \frac{\partial^2 U_m}{\partial z^2}+ (k^2n_0^2-\lambda) U_m = F_m.
\end{equation}
The solution to \eqref{8.12} which satisfies \eqref{8.5} is easily found,
$$
  U_m(\lambda, z)=\int\limits_{-\infty}^{+\infty}\!
  \frac{e^{i |z-\zeta|\sqrt{k^2n_0^2-\lambda}}}{2i\sqrt{k^2n_0^2-\lambda}}F_m(\lambda,\zeta)\,d\zeta,
$$
or if we use \eqref{8.11},
$$
  U_m(\lambda, z)=\int\limits_{-\infty}^{+\infty}\!\int\limits_{0}^{\infty}\!
  \frac{e^{i |z-\zeta|\sqrt{k^2n_0^2-\lambda}}}{2i\sqrt{k^2n_0^2-\lambda}}
  j_m(\rho,\lambda)\sqrt{\rho}f_m(\rho,\zeta)\,d\rho\,d\zeta.
$$
$U_m(\lambda, z)$ was defined by \eqref{8.4}. Using the inversion
formula \eqref{7.3} given in theorem \ref{th7.1} we can recover
$u_m(r,z),$
$$
  \sqrt{r}u_{m}(r,z)=\frac{1}{\pi} \int\limits_{-\infty}^{\infty}\!
  j_m(r,\lambda) U_m(\lambda, z)\, d\chi_m(\lambda),
$$
or
$$
  u_m(r,z)=\frac{1}{\pi\sqrt{r\rho}}
  \int\limits_{-\infty}^{\infty}\! \int\limits_{-\infty}^{+\infty}\!\int\limits_{0}^{\infty}\!
  \frac{e^{i |z-\zeta|\sqrt{k^2n_0^2-\lambda}}}{2i\sqrt{k^2n_0^2-\lambda}}
  j_m(\rho,\lambda) j_m(r,\lambda)\rho f_m(\rho,\zeta)
  \,d\rho\,d\zeta\, d\chi_m(\lambda).
$$
Now, to get \eqref{8.6} with $G(r, \rho; \vartheta, t; z, \zeta)$ given by
\eqref{8.7} we need to substitute $f_m(\rho,z)$ from \eqref{8.9}, find
$u(r, \vartheta, z)$ from \eqref{8.1} and interchange the order of
integration so that the inner-most integral is the one in respect to
$\lambda.$
\end{proof}

The theorem we just proved shows that the electromagnetic field generated by
the source $f$ can be decomposed in two parts: the guided part, which
is a sum of guided modes decaying in $r$ either exponentially or as a
power of $r,$ and a radiation part, which is obtained by summing in
$m$ and integrating in $\lambda.$ For each $m\in\mathbb Z$ the set of
guided modes is finite, as it was shown in lemma \ref{lm6.2}. The next
theorem will prove a stronger result.

\begin{theorem}\label{th8.3}
  The total number of guided modes (in all $m\in\mathbb Z$) is finite.
\end{theorem}
\begin{proof}
We just need to show that for $|m|$ large enough there are no more
guided modes.  A mode $j_m(r,\lambda)$ is guided, if
$$
  M_0^m(\lambda)-M_\infty^m(\lambda)=0.
$$
All $\lambda$ for which this equality happens are in $(0, d^2],$ as proved in
lemma \ref{lm6.2}. By using \eqref{6.1}, \eqref{6.7a} and \eqref{6.8} we can write the this equality
as
$$
  \frac{j_m'(R, \lambda)}{j_m(R, \lambda)}=\frac{k_m'(R, \lambda)}{k_m(R, \lambda)}, \mbox{ if } \lambda< d^2,
$$
or
$$
  \frac{j_m'(R, \lambda)}{j_m(R, \lambda)}=-\frac{|m|-1/2}{R}, \mbox{ if } \lambda=d^2,
$$
where $k(r,\lambda)$ is given by \eqref{4.5}.

The left-hand side of these equalities is strictly positive for
$|m|$ large, as it follows from lemma \ref{lm3.3}. Thus, the second of these
equalities is not possible.  To show that the first one cannot happen,
it suffices to prove that $k_m'(R, \lambda)< 0$ for $|m|\ge 1,$ since we know
that $k_m(r, \lambda)>0,$ for all $r>0.$

The function $r\to k_m(r, \lambda)$ will satisfy the equation
$$
  k_m''(r,\lambda)=\bigg\{d^2-\lambda + \frac{m^2-1/4}{r^2}\bigg\}k_m(r,\lambda),
$$
which implies that $k_m''(r,\lambda)>0$ for $r>0,$ and so,
$k_m'(r,\lambda)$ is an increasing function of $r.$ From this and
\eqref{4.8} it follows that $k_m'(r,\lambda)<0$ for all $r>0.$
This finishes the proof of the theorem.
\end{proof}

\section{Conclusion} In this paper we have constructed a framework
for analyzing waveguide problems which is based on a transform
theory. The construction of the transform was more difficult, but the
final form relatively similar, to the 2-D case.\cite{MS} The primary
tool in obtaining the transform was the theory of self-adjoint singular
eigenvalue problems.

This paper completes the study of the wave propagation in a infinite cylindrical waveguide. We obtained a Green's
function valid for every choice of the index of refraction of the core with cylindrical symmetry. In particular,
it is enough to solve \eqref{6.17a} for $\lambda\in (0,d^2]$ and the differential equation \eqref{2.5} in the core
region in order to obtain the corresponding Green's function.

The obtained formula for Green's functions is very amenable to
computation. In a future article we will calculate explicitly Green's
function in the cases of a \emph{step-index fiber} and a \emph{coaxial
waveguide} and will display numerical results.

\nonumsection{Acknowledgments}
We wish to express our gratitude to our advisers, Fadil Santosa
and Rolando Magnanini\footnote{Fadil Santosa is the adviser of O.
A. and Rolando Magnanini is the adviser of G. C.}, for proposing
this topic to us and for the many discussions and countless
suggestions while writing this paper.

\nonumsection{References}


\begin{thebibliography}{Z-Z-Z}
\bibitem[1]{AS}
   M. Abramowitz and I. A. Stegun,
  {\small\bf Handbook of mathematical functions} (Dover, New York, 1965).
\bibitem[2]{CL}
   E. A. Coddington and N. Levinson,
  {\small\bf Theory of Ordinary Differential Equations} (McGraw-Hill, New
  York, 1955).
\bibitem[3]{MS} 
  R. Magnanini and F. Santosa,
  {\small\it Wave propagation in a 2-D optical waveguide}, {\small\it SIAM J. Appl.
  Math.,} {\small\bf 61} (2001) 1237 -- 1252.
\bibitem[4]{mos} 
  W. Magnus, F. Oberhettiger and R. P. Soni,
  {\small\bf Formulas and theorems for the Special Functions of Mathematical Physics}
  (Springer-Verlag, Berlin, 1966).
\bibitem[5]{Ma} 
  D. Marcuse,
  {\small\bf Light Transmission Optics} (Van Nostrand Reinhold Company, New York, 1982).
\bibitem[6]{SL} 
  A. W. Snyder and D. Love,
  {\small\bf Optical Waveguide Theory} (Chapman and Hall, London, 1974).
\bibitem[7]{Ti} 
  E. C. Titchmarsh,
  {\small\bf Eigenfunction expansions associated with second-order differential equations}
  (Oxford at the Clarendon Press, Oxford, 1946).
\end{thebibliography}
\end{document}